\documentclass[11pt]{amsart}

\usepackage{epsfig}
\usepackage[]{rotating}

\usepackage[]{rotating}

\newtheorem{thm}{Theorem}[section]
\newtheorem{cor}[thm]{Corollary}
\newtheorem{lem}[thm]{Lemma}
\newtheorem{prop}[thm]{Proposition}

\theoremstyle{remark}
\newtheorem{rem}[thm]{Remark}

\newtheorem{nota}[thm]{Notation}

\newcommand{\N}{{\mathbb{N}}}

\newcommand{\sign}{{\textrm{sign}}}
\newcommand{\grad}{{\textrm{grad }}}
\newcommand{\ric}{{\textrm{Ric}}}
\newcommand{\bcwp}{$(\psi,\mu)$-\emph{bcwp}}
\newcommand{\bcwpbracket}[2]{$[#1,#2]$-\emph{bcwp}}

\newcommand{\bcwpar}[1]{$(#1)$-\emph{bcwp}}

\newcommand{\discr}[1]{\emph{discr }(#1)}
\newcommand{\twistpar}[4]{$#1 \times_{[#3,#4]} #2$}

\newcommand{\meas}{{\textrm{meas}}}

\numberwithin{equation}{section}
\renewcommand{\theequation}{\arabic{section}.\arabic{equation}}

\begin{document}

\title{About Curvature, Conformal Metrics and Warped Products}

\author{Fernando Dobarro \\
\& \\
B\"{u}lent \"{U}nal}

\address[F. Dobarro]{Dipartimento di Matematica e Informatica,
Universit\`{a} degli Studi di Trieste, Via Valerio 12/b, I-34127
Trieste, Italy} \email {dobarro@dmi.units.it}

\address[B. \"{U}nal]{Department of Mathematics, Bilkent University,
         Bilkent, 06800 Ankara, Turkey}
\email{bulentunal@mail.com}



\date{\today}



%
%

\maketitle

\begin{abstract} We consider the curvature of a family
of warped products of two pseduo-Riemannian manifolds $(B,g_B)$
and $(F,g_F)$ furnished with metrics of the form $c^{2}g_B \oplus
w^2 g_F$ and, in particular, of the type $w^{2 \mu}g_B \oplus w^2
g_F$, where $c, w \colon B \to (0,\infty)$ are smooth functions
and $\mu$ is a real parameter. We obtain suitable expressions for
the Ricci tensor and scalar curvature of such products that allow
us to establish results about the existence of Einstein or
constant scalar curvature structures in these categories. If
$(B,g_B)$ is Riemannian, the latter question involves nonlinear
elliptic partial differential equations with concave-convex
nonlinearities and singular partial differential equations of the
Lichnerowicz-York type among others.
\end{abstract}


\section{Introduction}
\label{sec:Introduction}

The main concern of this paper is the curvature of a special family
of warped pseudo-metrics on product manifolds. We introduce a
suitable form for the relations among the involved curvatures in
such metrics and apply them to the existence and/or construction of
Einstein and constant scalar curvature metrics in this family.

\medskip

Let $B=(B_m,g_{B})$ and $F=(F_k,g_{F})$ be two pseudo-Riemannian
manifolds of dimensions $m \geq 1$ and $k \geq 0,$ respectively and
also let $B \times F$ be the usual product manifold of $B$ and $F$.
For a given smooth function $w \in C^\infty_{>0}(B)=\{v \in
C^\infty(B): v(x)>0, \, \forall \,x \in B\}$, the \emph{warped
product} $B \times _w F = ({(B \times_w F)}_{m+k},g=g_{B} +
w^2g_{F})$ was defined by Bishop and O'Neill in \cite{BishopONeil69}
in order to study manifolds of negative curvature.

\medskip

In this article, we deal with a particular class of warped products,
i.e. when the pseudo-metric in the base is affected by a conformal
change. Precisely, for given smooth functions $c, w \in
C^\infty_{>0}(B)$ we will call $({(B \times F)}_{m+k},g=c^2 g_{B} +
w^2g_{F})$ as a $[c,w]$-base conformal warped product (briefly
\bcwpbracket{c}{w}), denoted by \twistpar{B}{F}{c}{w}. We will
concentrate our attention on a special subclass of this structure,
namely when there is a relation between the \emph{conformal factor}
$c$ and the \emph{warping function} $w$ of the form $c=w^\mu$, where
$\mu $ is a real parameter and we will call the
\bcwpbracket{\psi^\mu}{\psi} as a \bcwpar{\psi,\mu}. Note that we
generically called the latter case as special base conformal warped
products, briefly \emph{sbcwp} in \cite{DobarroUnal04-2}.

\medskip

As we will explain in \S \ref{sec:main results}, metrics of this
type play a relevant role in several topics of differential
geometry and theoretical physics (see also
\cite{DobarroUnal04-2}). This article concerns curvature related
questions of these metrics which are of interest
not only in the applications, but also from the points of view of
differential geometry and the type of the involved nonlinear
partial differential equations (PDE), such as those with
concave-convex nonlinearities and the Lichnerowicz-York equations.

\medskip

The article is organized in the following way: in \S \ref{sec:main
results} after a brief description of several fields where
pseudo-metrics described as above are applied, we formulate the
curvature problems that we deal within the next sections and give
the statements of the main results. In \S \ref{sec:The curvature
relations}, we state Theorems \ref{thm:ricci sbcwp m ge 3} and
\ref{thm:scurv conf warped m ge 2} in order to express the Ricci
tensor and scalar curvature of a \bcwpar{\psi,\mu} and sketch their
proofs (see \cite[Section 3]{DobarroUnal04-2} for detailed
computations). In \S \ref{sec:The problem (Pb-sc)} and
\ref{sec:concave-convex}, we establish our main results about the
existence of \bcwpar{\psi,\mu}'s of constant scalar curvature with
compact Riemannian base.


\section{Motivations and Main results}
\label{sec:main results}

As we announced in the introduction, we firstly want to mention some
of the major fields of differential geometry and theoretical physics
where base conformal warped products are applied.

\begin{description}
\item[i] In the construction of a large class of non trivial
static anti de Sitter vacuum space-times
\begin{itemize}
  \item
In the Schwarzschild solutions of the Einstein equations (see
\cite{Anderson-Chrusciel-Delay02,Besse87,
Hawking-Ellis73,ONeil83,Schwarzschild1916,Thorne03}).
  \item
In the Riemannian Schwarzschild metric, namely (see
\cite{Anderson-Chrusciel-Delay02}).
  \item
In the \emph{``generalized Riemannian anti de Sitter
$\mathbf{T}^{2}$ black hole metrics"} (see \S3.2 of
\cite{Anderson-Chrusciel-Delay02} for details).
    \item
In the Ba\~{n}ados-Teitelboim-Zanelli (BTZ) and de Sitter (dS)
black holes (see \cite{Aharony-Gubser-Maldacena-Ooguri-Oz,
Banados-Teitelboim-Zanelli92,
Banados-Henneaux-Teitelboim-Zanelli93,
DobarroUnal04,Hong-Choi-Park03} for details).
\end{itemize}

Indeed, all of them can be generated by an approach of the
following type: let $(F_2,g_F)$ be a pseudo-Riemannian manifold
and $g$ be a pseudo-metric on $\mathbb{R}_+ \times \mathbb{R}
\times F_2$ defined by
\begin{equation}\label{eq:Schwarzschild to bcwp-1}
g= \frac{1}{u^{2}(r)}dr^{2} \pm u^{2}(r) dt^2 + r^{2}g_F.
\end{equation}
After the change of variables $s=r^2$, $\displaystyle
y=\frac{1}{2} t$, there results $ds^2=4 r^2 dr^2$ and $
\displaystyle dy^2 = \frac{1}{4} dt^2$. Then
\eqref{eq:Schwarzschild to bcwp-1} is equivalent to
\begin{equation}
\label{eq:Schwarzschild to bcwp-2}
\begin{split}
g&= \frac{1}{\sqrt{s}}\Big[\frac{1}{4 \sqrt{s}
u^{2}(\sqrt{s})}ds^{2} \pm 4 \sqrt{s} u^{2}(\sqrt{s}) dy^2\Big] +
s g_F \\
&= (s^\frac{1}{2})^{2(-\frac{1}{2})}
\Big[(2 s^\frac{1}{4}  u(s^\frac{1}{2}))^{2(-1)}ds^{2} \pm
(2 s^\frac{1}{4}  u(s^\frac{1}{2}))^2 dy^2\Big]
 + (s^\frac{1}{2})^{2} g_F.
\end{split}
\end{equation}
Note that roughly speaking, $g$ is a nested application of two
\bcwp's. That is, on $\mathbb{R}_+ \times \mathbb{R}$ and taking
\begin{equation}\label{eq:Schwarzschild to bcwp-3}
  \psi_1(s) = 2 s^\frac{1}{4}  u(s^\frac{1}{2}) \textrm{ and }
\mu_1=-1,
\end{equation}
the metric inside the brackets in the last member of
\eqref{eq:Schwarzschild to bcwp-2} is a \bcwpar{\psi_1,\mu_1},
while the metric $g$ on $(\mathbb{R}_+ \times \mathbb{R}) \times
F_2$ is a \bcwpar{\psi_2,\mu_2} with
\begin{equation}\label{eq:Schwarzschild to bcwp-4}
  \psi_2(s,y)=s^\frac{1}{2} \textrm{ and } \mu_2 = -\frac{1}{2}.
\end{equation}
In the last section of \cite{DobarroUnal04-2}, through the
application of Theorems \ref{thm:ricci sbcwp m ge 3} and
\ref{thm:scurv conf warped m ge 2} below and several standard
computations, we generalized the latter approach to the case of an
Einstein fiber $(F_k,g_F)$ with dimension $k \ge 2$.

\item [ii] In the study of the equivariant isometric
embeddings of space-time slices in Minkowski spaces (see
\cite{Giblin-Marlof-Garvey03, Giblin-Hwang04}).

\item[iii] In the Kaluza-Klein theory (see \cite[\S7.6,
Particle Physics and Geometry]{Wesson99}, \cite{OverduinWesson98}
and \cite{Wesson01}) and in the Randall-Sundrum theory
\cite{Frolov01,Greene-Schalm-Shiu00,Randall-Sundrum99a,Randall-Sundrum99b,Randjbar
Daemi-Rubakov04,Soda02} with $\mu$ as a free parameter. For example,
in \cite{Ito01} the following metric is considered
\begin{equation}\label{eq:Ito1}
  e^{2\mathcal{A}(y)} g_{ij}dx^{i}dx^{j} + e^{2\mathcal{B}(y)}dy^{2},
\end{equation}
with the notation $\{x^{i}\}$, $i= 0,1,2,3$ for the coordinates in
the 4-dimensional space-time and $x^{5}=y$ for the fifth
coordinate on an extra dimension. In particular, Ito takes the
ansatz
\begin{equation}\label{eq:Ito2}
  \mathcal{B}=\alpha \mathcal{A},
\end{equation}
which corresponds exactly to our \emph{sbcwp} metrics,
considering $g_{B}=dy^{2}$, $g_{F}=g_{ij}dx^{i}dx^{j}$,
$\displaystyle \psi(y) = e^{\frac{\mathcal{B}(y)}{\alpha}} =
e^{\mathcal{A}(y)}$ and $\displaystyle \mu = \alpha$.

\item [iii] In String and Supergravity theories, for instance, in
the Maldacena conjecture about the duality between compactifications
of M/string theory on various Anti-de Sitter space-times and various
conformal field theories (see \cite{Maldacena98,Petersen99}) and in
warped compactifications (see \cite{Greene-Schalm-Shiu00,
Strominger86} and references therein). Besides all of these, there
are also frequent occurrences of this type of metrics in string
topics (see \cite{
Gauntlett-Kim-Waldram01, Gauntlett-Kim-Waldram01-b,
Gauntlett-Kim-Pakis-Waldram02,
Gauntlett-Martelli-Sparks-Waldram04,GhezelbashMann04, Lidsey00,
PapadopoulosTownsend99,Soda02} and also
\cite{Aharony-Gubser-Maldacena-Ooguri-Oz, Argurio98, Schmidt04}
for some reviews about these topics).

\item [iv] In the derivation of effective theories for warped
compactification of supergravity and the Ho\u{r}ava-Witten model
(see \cite{Kodama-Uzawa05, Kodama-Uzawa06}). For instance, in
\cite{Kodama-Uzawa06} the ansatz $ds^2=h^\alpha ds^2(X_4) +
h^\beta ds^2(Y)$ is considered where $X_4$ is a four-dimensional
space-time with coordinates $x^\mu$, $Y$ is a Calabi-Yau manifold
(the so called internal space) and $h$ depends on the
four-dimensional coordinates $x^\mu$, in order to study the
dynamics of the four-dimensional effective theory. We note that in
those articles, the structure of the expressions of the Ricci
tensor and scalar curvature of the involved metrics result
particularly useful. We observe that they correspond to very
particular cases of the expressions obtained by us in
\cite{DobarroUnal04-2}, see also Theorems \ref{thm:ricci sbcwp m
ge 3} and \ref{thm:scurv conf warped m ge 2} and Proposition
\ref{prop:Einstein condition} stated below.

\item [v] In the discussion of Birkhoff-type theorems (generally
speaking these are the theorems in which the gravitational vacuum
solutions admit more symmetry than the inserted metric ansatz, (see
\cite[page 372]{Hawking-Ellis73} and  \cite[Chapter
3]{Beem-Ehrlich-Easley96}) for rigorous statements), especially in
Equation 6.1 of \cite{Schmidt97} where, H-J. Schmidt considers a
special form of a \emph{bcwp} and basically shows that if a
\emph{bcwp} of this form is Einstein, then it admits one Killing
vector more than the fiber. In order to achieve that, the author
considers for a specific value of $\mu$, namely $\mu = (1-k)/2$, in
the following problem:

\begin{quote}
\emph{Does there exist a smooth function $\psi \in
C^\infty_{>0}(B)$ such that the corresponding \bcwpar{\psi,\mu}
$(B_2 \times F_k,\psi^{2 \mu}g_{B} + \psi^{2} g_{F})$ is an
Einstein manifold?} (see also \textbf{(Pb-Eins.)} below.)
\end{quote}

\item [vi] In the study of bi-conformal transformations,
bi-conformal vector fields and their applications (see
\cite[Remark in Section 7]{GarciaParrado-Senovilla04} and
\cite[Sections 7 and 8]{GarciaParrado04}).

\item[vii] In the study of the spectrum of the Laplace-Beltrami
operator for $p-$forms. For instance in Equation $(1.1)$ of
\cite{Antoci03}, the author considers the structure that follows:
let $\overline{M}$ be an $n$-dimensional compact, Riemannian
manifold with boundary, and let $y$ be a boundary-defining
function; she endows the interior $M$ of $\overline{M}$ with a
Riemannian metric $ds^2$ such that in a small tubular neighborhood
of $\partial M$ in $M$, $ds^2$ takes the form
\begin{equation}\label{eq:antoci}
ds^2 = e^{-2(a+1)t} dt^2 + e^{-2bt}d\theta^2_{\partial M},
\end{equation}
where $t := -\log y \in (c,+\infty)$ and $d\theta^2_{\partial M}$
is the Riemannian metric on $\partial M$ (see \cite{Antoci03,
Melrose95} and references therein for details).
\end{description}

\medskip

\begin{nota}
From now on, we will use the Einstein summation convention over
repeated indices and consider only connected manifolds.
Furthermore, we will denote the Laplace-Beltrami operator on a
pseudo-Riemannian manifold  $(N,h)$ by $\Delta_{N}(\cdot),$ i.e.,
$\Delta_{N}(\cdot)={\nabla^{N}}^i{\nabla^{N}}_i (\cdot).$ Note
that $\Delta_{N}$ is elliptic if $(N,h)$ is Riemannian and it is
hyperbolic when $(N,h)$ is Lorentzian. If $(N,h)$ is neither
Riemannian nor Lorentzian, then the operator is ultra-hyperbolic.

\medskip

\noindent Furthermore, we will consider the Hessian of a function $v
\in C^\infty(N)$, denoted by $H_{h}^{v}$ or $H_{N}^{v}$, so that the
second covariant differential of $v$ is given by
$H_{h}^{v}=\nabla(\nabla v)$. Recall that the Hessian is a symmetric
$(0,2)$ tensor field satisfying
\begin{equation}
 \label{eq:hessian.0}
 H_{h}^{v} (X,Y) = X Y v - (\nabla_{X}Y)v = h(\nabla_{X}(\grad
 v),Y),
\end{equation}
for any smooth vector fields $X,Y$ on $N.$

\medskip

\noindent For a given pseudo-Riemannian manifold $N=(N,h)$ we will
denote its Riemann curvature tensor, Ricci tensor and scalar
curvature by $R_N$, $\ric_N$ and $S_N$, respectively.

\medskip

\noindent We will denote the set of all lifts of all vector fields
of $B$ by $\mathfrak L(B).$ Note that the lift of a vector field
$X$ on $B$ denoted by $\widetilde X$ is the vector field on $B
\times F$ given by ${\rm d} \pi (\widetilde{X})=X$ where $\pi
\colon B \times F \to B$ is the usual projection map.
\end{nota}

\medskip

In Section \ref{sec:The curvature relations}, we will sketch the
proofs of the following two theorems related to the Ricci tensor and
the scalar curvature of a generic \bcwpar{\psi,\mu}.

\begin{thm}
\label{thm:ricci sbcwp m ge 3} Let $B = (B_m,g_{B})$ and $F =
(F_k,g_{F})$ be two pseudo-Rieman\-nian manifolds with $m \geq 3$
and $k \geq 1$, respectively and also let $\displaystyle \mu \in
\mathbb{R}\setminus \{0, 1, \overline{\mu},
\overline{\mu}_{\pm}\}$ be a real number with
\begin{equation*}\label{eq: Ricci special values}
    \overline{\mu}:=\displaystyle-\frac{k}{m-2} \textit{ and }
    \overline{\mu}_{\pm}:=\overline{\mu}\pm\sqrt{\overline{\mu}^{2}-\overline{\mu}}.
\end{equation*}
%
Suppose $\psi \in C^\infty_{>0}(B)$. Then the Ricci curvature
tensor of the corresponding \emph{\bcwpar{\psi,\mu}}, denoted by
\ric ~ verifies the relation
\begin{equation}\label{eq:ricci special 1}
\begin{split}
&\displaystyle{{\ric}={\ric}_B+
\beta^{H}\frac{1}{\psi^{\frac{1}{\alpha ^{H}}}}{\rm
H}_{B}^{\psi^{\frac{1}{\alpha ^{H}}}} - \beta^{\Delta}
\frac{1}{\psi ^{\frac{1}{\alpha ^{\Delta}}}} {\Delta}_{B}
\psi^{\frac{1}{\alpha ^{\Delta}}} }g_B
\textrm{ on } \mathcal{L}(B)\times \mathcal{L}(B),
\\
&\displaystyle{{\ric}=0} \textrm{ on } \mathcal{L}(B)\times
\mathcal{L}(F),
\\
&\displaystyle{{\ric}={\ric}_F
-\frac{1}{ \psi^{2(\mu - 1)}}
\frac{\beta^{\Delta}}{\mu} \frac{1}{\psi ^{\frac{1}{\alpha
^{\Delta}}}} {\Delta}_{B} \psi^{\frac{1}{\alpha ^{\Delta}}}
g_F}
\textrm{ on } \mathcal{L}(F)\times \mathcal{L}(F),
\end{split}
\end{equation}
where
\begin{equation}\label{eq:alphaDH betaDH particular}
  \begin{array}{lcl}
    \alpha^{\Delta} &=&   \displaystyle \frac{1}{(m-2)\mu + k},      \\
    \beta^{\Delta} &=& \displaystyle \frac{\mu}{(m-2)\mu + k},      \\
    \alpha^{H} &=& \displaystyle \frac{-[(m-2)\mu + k]}{\mu [(m-2)\mu
+ k]+k(\mu - 1)},    \\
    \beta^{H} &=&  \displaystyle \frac{[(m-2)\mu + k]^{2}}{\mu
[(m-2)\mu + k]+k(\mu - 1)}.
  \end{array}
\end{equation}

\end{thm}

\medskip

\begin{thm} \label{thm:scurv conf warped m ge 2}
Let $B = (B_m,g_{B})$ and $F = (F_k,g_{F})$ be two
pseudo-Rieman\-nian manifolds of dimensions $m \geq 2$ and $k \geq
0,$ respectively. Suppose that $S_B$ and $S_F$ denote the scalar
curvatures of $B = (B_m,g_{B})$ and $F = (F_k,g_{F}),$
respectively. If $\displaystyle \mu \in \mathbb{R}$ and $\psi \in
C^\infty_{>0}(B)$, then the scalar curvature $S$ of the
corresponding \emph{\bcwpar{\psi,\mu}} verifies,
\begin{itemize}
  \item[\bf(i)] If $\displaystyle \mu \neq -\frac{k}{m-1}$, then
\begin{equation}
-\beta \Delta_{B}u + S_{B} u = S u^{2 \mu \alpha + 1}   -
S_{F}u^{2(\mu - 1)\alpha + 1}
 \label{eq:sc principal 1}
\end{equation}
where
\begin{equation}
 \label{eq:alpha principal}
 \alpha =
\frac{2[k+(m-1)\mu]}{\{[k+(m-1)\mu]+(1-\mu)\}k + (m-2)\mu
[k+(m-1)\mu]},
\end{equation}
\begin{equation}
 \label{eq:beta principal}
  \beta = \alpha 2[k+(m-1)\mu]>0
\end{equation}
and $\psi = u^{\alpha}>0$.
   \item[\bf(ii)] If $\displaystyle \mu = -\frac{k}{m-1}$, then
\begin{equation}
\label{eq:sc special}
%
%
  -k\left[1 + \frac{k}{m-1}\right]
  \frac{|\nabla^{B}
  \psi|_{B}^{2}}{\psi^{2}}=
%
%
  \psi^{-2\frac{k}{m-1}} [S  -    S_{F} \psi^{-2}]- S_{B}.
\end{equation}
\end{itemize}
\end{thm}

%
From the mathematical and physical points of view, there are
several interesting questions about \bcwpar{\psi,\mu}'s. In
\cite{DobarroUnal04-2} we began the study of existence and/or
construction of Einstein \bcwpar{\psi,\mu}'s and those of constant
scalar curvature. These questions are closely connected to
Theorems \ref{thm:ricci sbcwp m ge 3} and \ref{thm:scurv conf
warped m ge 2}.

\medskip

In \cite{DobarroUnal04-2}, by applying Theorem \ref{thm:ricci
sbcwp m ge 3}, we give suitable conditions that allow us to study
some particular cases of the problem:

\medskip
\begin{quote}
  \textbf{(Pb-Eins.)} Given $\mu \in \mathbb{R}$, does there
  exist a smooth function $\psi \in C^\infty_{>0}(B)$ such that
  the corresponding
\bcwpar{\psi,\mu} is an Einstein manifold?
\end{quote}
\medskip
In particular, we obtain the following result as an immediate
corollary of Theorem \ref{thm:ricci sbcwp m ge 3}.
\begin{prop}\label{prop:Einstein condition}
Let us assume the hypothesis of Theorem \ref{thm:ricci sbcwp m ge
3}. Then the corresponding \emph{\bcwpar{\psi,\mu}} is an Einstein
manifold with $\lambda$ if and only if $(F,g_F)$ is Einstein with
$\nu$ constant and the system that follows is verified
\begin{equation}\label{eq:einstein: special m=m}
\begin{split}
&\displaystyle{\lambda \psi^{2\mu}g_B={\rm Ric}_B+
\beta^{H}\frac{1}{\psi^{\frac{1}{\alpha ^{H}}}}{\rm
H}_{B}^{\psi^{\frac{1}{\alpha ^{H}}}} - \beta^{\Delta}
\frac{1}{\psi ^{\frac{1}{\alpha ^{\Delta}}}} {\Delta}_{B}
\psi^{\frac{1}{\alpha ^{\Delta}}} }g_B
\textrm{ on } \mathcal{L}(B)\times \mathcal{L}(B)
%
%
\\
&\displaystyle{\lambda \psi^2=\nu
-\frac{1}{ \psi^{2(\mu - 1)}}
\frac{\beta^{\Delta}}{\mu} \frac{1}{\psi ^{\frac{1}{\alpha
^{\Delta}}}} {\Delta}_{B} \psi^{\frac{1}{\alpha ^{\Delta}}}
},
%
%
\end{split}
\end{equation}
where the coefficients are given by \eqref{eq:alphaDH betaDH
particular}.
\end{prop}

\noindent Compare the system \eqref{eq:einstein: special m=m} with
the well known one for a classical warped product in \cite{Besse87,
Kim-Kim01, ONeil83}. By studying \eqref{eq:einstein: special m=m},
we have obtained the generalization of the construction exposed in
the above motivational examples in \textbf{i} and \textbf{v}, among
other related results. We suggest the interested reader consider the
results about the problem \textbf{(Pb-Eins.)} stated in
\cite{DobarroUnal04-2}.

\medskip

Now, we focus on the problems which we will deal in \S \ref{sec:The
problem (Pb-sc)}. Let $B=(B_m,g_B)$ and $F=(F_k,g_F)$ be
pseudo-Riemannian manifolds.

\noindent There is an extensive number of publications about the
well known Yamabe problem namely:

\medskip
\begin{quote}
\label{(Ya):problem}
    \textbf{(Ya)} \cite{Yamabe63,Trudinger68,Schoen84,Aubin98}
    Does there exist a function
    $\varphi \in C^\infty_{>0}(B)$ such that
    $(B_m,\varphi^{\frac{4}{m-2}}g_{B})$ has
    constant scalar curvature?
\end{quote}
\medskip

\noindent Analogously, in several articles the following problem
has been studied:

\medskip
\begin{quote}
\label{(cscwp):problem}
    \textbf{(cscwp)} \cite{DobarroLamiDozo87}
    Is there a function $w \in C^\infty_{>0}(B)$
    such that the warped product
    $B \times_{w} F$ 
    has constant scalar curvature?
\end{quote}

\medskip

\noindent \textit{In the sequel we will suppose that $B=(B_m,g_B)$
is a Riemannian manifold.}

\medskip

\noindent Thus, both problems bring to the study of the existence
of positive solutions for nonlinear elliptic equations on
Riemannian manifolds. The involved nonlinearities are powers with
Sobolev critical exponent for the Yamabe problem and sub-linear
(linear if the dimension $k$ of the fiber is $3$) for the problem
of constant scalar curvature of a warped product.

\medskip

\noindent In Section \ref{sec:The problem (Pb-sc)},
we deal with a mixed problem between \textbf{(Ya)} and
\textbf{(cscwp)} which is already proposed in
\cite{DobarroUnal04-2}, namely:

\medskip
\begin{quote}
  \textbf{(Pb-sc)} Given $\mu \in \mathbb{R}$, does there
  exist a function $\psi \in C^\infty_{>0}(B)$ such that
  the corresponding \bcwpar{\psi,\mu}
  has constant scalar curvature?
\end{quote}

\medskip


\noindent \textsl{Note that when $\mu = 0,$ \textbf{(Pb-sc)}
corresponds to the problem \textbf{(cscwp)}, whereas when the
dimension of the fiber $k=0$ and $\mu = 1$, then \textbf{(Pb-sc)}
corresponds to \textbf{(Ya)} for the base manifold. Finally
\textbf{(Pb-sc)} corresponds to \textbf{(Ya)} for the usual
product metric with a conformal factor in $C^\infty_{>0}(B)$ when
$\mu=1$.}

\medskip

\noindent Under the hypothesis of \textbf{Theorem \ref{thm:scurv
conf warped m ge 2} i}, the analysis of the problem
\textbf{(Pb-sc)} brings to the study of the existence and
multiplicity of solutions $u \in C^\infty_{>0}(B)$ of
\begin{equation}
   -\beta \Delta_{B}u + S_{B} u =
  \lambda u^{2 \mu \alpha + 1} - S_{F}u^{2(\mu - 1)\alpha + 1},
\label{eq:scalar curv bcwp1}
\end{equation}
where all the components of the equation are like in \textbf{Theorem
\ref{thm:scurv conf warped m ge 2} i} and $\lambda$ (the conjectured
constant scalar curvature of the corresponding \bcwpar{\psi,\mu}) is
a real parameter. We observe that an easy argument of separation of
variables, like in \cite[\S2]{CotiZelatiDobarroMusina97} and
\cite{DobarroLamiDozo87}, shows that there exists a positive
solution of \eqref{eq:scalar curv bcwp1} only if the scalar
curvature of the fiber $S_H$ is constant. Thus this will be a
natural assumption in the study of \textbf{(Pb-sc)}.

\noindent Furthermore, note that the involved nonlinearities in the
right hand side of \eqref{eq:scalar curv bcwp1} dramatically change
with the choice of the parameters, an exhaustive analysis of these
changes is the subject matter of \cite[\S6]{DobarroUnal04-2}.

\noindent There are several partial results about semi-linear
elliptic equations like \eqref{eq:scalar curv bcwp1} with
different boundary conditions, see for instance
\cite{Alama99,Ambrosetti-Brezis-Cerami94,
Ambrosetti-GarciaAzorero-Peral00,Ambrosetti-Rabinowitz73,
Chabrowski-doO02,CortazarElguetaFelmer93,De
Figueiredo-Gossez-Ubilla03,Taylor96,Willem96} and references in
\cite{DobarroUnal04-2}.

\smallskip

\noindent \emph{In this article we will state our first results
about the problem \textbf{(Pb-sc)} when the base $B$ is a compact
Riemannian manifold of dimension $m \ge 3$ and the fiber $F$ has
non-positive constant scalar curvature $S_F$.}

\smallskip

\noindent For brevity of our study, it will be useful to introduce
the following notation: $\displaystyle{\mu_{sc}:=\mu_{sc}(m,k)= }$
$\displaystyle{ -\frac{k}{m-1}}$ and $\mu_{p_Y}=\mu_{p_Y}(m,k):=
\displaystyle{-\frac{k+1}{m-2}}$ ($sc$ as scalar curvature and $Y$
as Yamabe). Notice that $\mu_{p_Y}<\mu_{sc}<0$.

\smallskip

\noindent \textit{We plan to study the case of $\mu=\mu_{sc}$ in a
preceding project, therefore the related results are not going to be
presented here.}

\smallskip

\noindent We can synthesize our results about $\textbf{(Pb-sc)}$
in the case of non-positive $S_F$ as follow.

\smallskip

\begin{itemize}
\item\emph{The case  of scalar flat fiber, i.e. $S_F=0$.}

\begin{thm} \label{thm:Pb-sc SF0} If $\mu \in
(\mu_{p_Y},\mu_{sc}) \, \cup \, (\mu_{sc},+\infty) $
the answer to
$\textbf{(Pb-sc)}$ is affirmative.
\end{thm}

\noindent By assuming some additional restrictions on the scalar
curvature of the base $S_B$, we obtain existence results for the
range $\mu \le \mu_{p_Y}$.
%
%
%
%
\end{itemize}

\smallskip

\begin{itemize}\item \emph{The case of fiber with negative constant scalar
curvature, i.e. $S_F<0$.}


\noindent In order to describe the $\mu-$ranges of validity of the
results, we will apply the notations introduced in
\cite[\S5]{DobarroUnal04-2} (see Appendix \ref{appendix:A} for a
brief introduction of these notations).

\begin{thm}\label{thm:Pb-sc SF<0} If ``$(m,k) \in D$ and $\mu \in (0,1)$" or
``$(m,k) \in \mathcal{C}D$ and $\mu \in (0,1) \cap (\mu_{-},
\mu_{+})$" or ``$(m,k) \in \mathcal{C}D$ and $\mu \in (0,1) \cap
\; \mathcal{C}[\mu_{-} , \mu_{+}]$", then the answer to
$\textbf{(Pb-sc)}$ is affirmative.
\end{thm}

\begin{rem}\label{rem:ABC-YL}
The first two cases in Theorem \ref{thm:Pb-sc SF<0} will be
studied by adapting the ideas in \cite{Ambrosetti-Brezis-Cerami94}
and the last case by applying the results in \cite[p.
99]{Taylor96}. In the former - Theorem \ref{thm:concave-convex},
the involved nonlinearities correspond to the so called
concave-convex whereas in the latter - Theorem \ref{thm:p superlin
q sing}, they are singular as in the Lichnerowicz-York equation
about the constraints for the Einstein equations (see
\cite{Choquet-Bruhat}, \cite{Hebey-Pacard-Pollack2007},
\cite{Murchadha}, \cite[p. 542-543]{GRA} and
\cite[Chp.18]{Taylor96}).

\noindent Similarly to the case of $S_F=0$, we obtain existence
results for some remaining $\mu-$ranges by assuming some
additional restrictions for the scalar curvature of the base
$S_B$.
\end{rem}
\end{itemize}

\noindent Naturally the study of $\textbf{(Pb-sc)}$ allows us to
obtain partial results of the related question:
\begin{quote}
  Given $\mu \in \mathbb{R}$ and $\lambda \in \mathbb{R}$ does there
  exist a function $\psi \in C^\infty_{>0}(B)$ such that
  the corresponding \bcwpar{\psi,\mu}
  has constant scalar curvature $\lambda$?
\end{quote}
These are stated in the several theorems and propositions in \S
\ref{sec:The problem (Pb-sc)}.

\section{The curvature relations - Sketch of the proofs}
\label{sec:The curvature relations}

The proofs of Theorems \ref{thm:ricci sbcwp m ge 3} and
\ref{thm:scurv conf warped m ge 2} require long and yet standard
computations of the Riemann and Ricci tensors and the scalar
curvature of a general base conformal warped product. Here, we
reproduce the results for the Ricci tensor and the scalar curvature,
and we also suggest the reader see \cite[\S3]{DobarroUnal04-2} for
the complete computations.

\begin{thm} \label{thm:global riki}
The Ricci tensor of \bcwpbracket{c}{w}, denoted by $Ric$ satisfies

\begin{enumerate}
\item $\displaystyle{{\ric}={\ric}_B-
\left[(m-2)\frac{1}{c}{\rm H}_{B}^c + k\frac{1}{w} {\rm H}_{B}^w
\right]}$\\
\mbox

\hspace{0.67 cm}  $\displaystyle{ +2(m-2)\frac{1}{c^{2}} dc
\otimes dc +k\frac{1}{wc} [dc \otimes dw + dw
\otimes dc]}$\\
\mbox

\hspace{0.67 cm} $\displaystyle -\left[(m-3)\frac{g_B(\nabla^B c
 ,\nabla^B c )}{c^{2}}+\frac{\Delta_{B}c}{c} + k
\frac{g_B(\nabla^B w ,\nabla^B c )}{wc}
\right]g_B $\\
\mbox

\hspace{0.67 cm} $\displaystyle \textrm{ on } \mathcal{L}(B)\times
\mathcal{L}(B)$, \mbox

\item $\displaystyle{{\ric}=0} \textrm{ on }
\mathcal{L}(B)\times \mathcal{L}(F)$,

\item $\displaystyle{{\ric}={\ric}_F
-\frac{w^{2}}{c^{2}} \left[ (m-2)
\frac{g_B(\nabla^B w ,\nabla^B c )}{wc} +
\frac{\Delta_{B}w}{w}\right.}$ \\
\mbox

\hspace{0.67 cm} $\displaystyle{\left.+ (k-1) \frac{g_B(\nabla^B w
,\nabla^B w )}{w^2}\right]  g_F} \textrm{ on }
\mathcal{L}(F)\times \mathcal{L}(F)$.
\end{enumerate}
\end{thm}

\begin{thm}
\label{thm:sca-c}
The scalar curvature $S$ of a \bcwpbracket{c}{w} is given by
\begin{eqnarray*}
c^{2}S & = & S_B + S_F\frac{c^{2}}{w^2}
-2(m-1)\frac{\Delta_{B}c}{c}- 2k
\frac{\Delta_{B} w}{w}\\
& - & (m-4)(m-1)   \frac{g_B(\nabla^B c ,\nabla^B c )}{c^{2}}\\
& - & 2k(m-2) \frac{g_B(\nabla^B w ,\nabla^B c )}{w c} \\
& - & k(k-1)\frac{g_B(\nabla^B w ,\nabla^B w )}{w^{2}}.
\end{eqnarray*}
\end{thm}

The following two lemmas (\ref{m-lem} and \ref{lem:hessian}) play a
central role in the proof of Theorems \ref{thm:ricci sbcwp m ge 3}
and \ref{thm:scurv conf warped m ge 2}. Indeed, it is sufficient to
apply them in a suitable mode and make use of Theorems
\ref{thm:global riki} and \ref{thm:sca-c} several times, the reader
can find all the details in \cite[\S 2 and 4]{DobarroUnal04-2}.

\medskip

Let $N = (N_n,h)$ be a pseudo-Rieman\-nian manifold of dimension
$n$,
%
$|\nabla (\cdot)|^{2}=|\nabla^{N} (\cdot)|_{N}^{2}=
h(\nabla^{N}(\cdot),\nabla^{N}(\cdot))$ and
$\Delta_{h}=\Delta_{N}$.

\begin{lem} \label{m-lem} Let $L_{h}$ be a differential operator
on $C^\infty_{>0}(N)$ defined by
\begin{equation} L_{h} v =\displaystyle\sum_{i=1}^k
r_{i}\frac{\Delta_{h}v^{a_{i}}}{v^{a_{i}}}, \label{m-lem1}
\end{equation}

where $r_{i},a_{i} \in \mathbb{R}$ and $\zeta :=\displaystyle
\sum_{i=1}^k r_{i}a_{i}$, $\eta:= \displaystyle \sum_{i=1}^k
r_{i}a_{i}^{2} $. Then,

\begin{itemize}
\item[\bf (i)] \begin{equation} L_{h} v = (\eta -
\zeta)\frac{\|\grad_h v \|_h^{2}}{v^{2}} + \zeta
\frac{\Delta_{h}v}{v}. \label{m-lem2}
\end{equation}

\item[\bf (ii)] If $\zeta \neq 0$ and $\eta \neq 0$, for $\alpha =
\displaystyle \frac{\zeta}{\eta}$ and $\beta = \displaystyle
\frac{\zeta^{2}}{\eta},$ then we have
\begin{equation} L_{h} v= \beta \frac
{\Delta_{h}v^{\frac{1}{\alpha}}}{v^{\frac{1}{\alpha}}}.
\label{m-lem3}
\end{equation}
\end{itemize}
\end{lem}

\begin{rem}\label{rem: curvature of multiply warped}
We also applied the latter lemma in the study of curvature of
multiply warped products (see  \cite{DobarroUnal04}).
\end{rem}

\begin{cor}
\label{cor:2.1}
    Let $L_{h}$ be a differential operator defined by
    \begin{equation}
     L_{h} v = r_{1}\frac{\Delta_{h}v^{a_{1}}}{v^{a_{1}}} +
           r_{2}\frac{\Delta_{h}v^{a_{2}}}{v^{a_{2}}}  \textrm{ for }  v \in
           C^\infty_{>0}(N),
        \label{eq:2.3}
    \end{equation}
    where $ r_{1}a_{1} +
r_{2}a_{2}\neq 0 $ and
    $ r_{1}a_{1}^{2} + r_{2}a_{2}^{2}\neq 0 $. Then, by changing
the variables $v = u^{\alpha}$ with $0 < u \in
    C^\infty(N)$,
    $\alpha = \displaystyle
    \frac{ r_{1}a_{1}+r_{2}a_{2}}
    { r_{1}a_{1}^{2} + r_{2}a_{2}^{2}}$ and
    $\beta = \displaystyle
    \frac{(r_{1}a_{1} + r_{2}a_{2})^{2}}
    {r_{1}a_{1}^{2} + r_{2}a_{2}^{2}}=\alpha (r_{1}a_{1} +
r_{2}a_{2})$ there results
    \begin{equation}
     L_{h} v= \beta \frac{\Delta_{h}u}{u}.
    \label{eq:2.4}
    \end{equation}
\end{cor}

\begin{rem}
\label{rem:2.1} By the change of variables as in Corollary
\ref{cor:2.1} equations of the type
\begin{equation}
L_{h} v = r_{1}\frac{\Delta_{h}v^{a_{1}}}{v^{a_{1}}} +
           r_{2}\frac{\Delta_{h}v^{a_{2}}}{v^{a_{2}}} = H(v,x,s),
    \label{eq:nonlinear eq}
\end{equation}
transform into
\begin{equation}
\beta \Delta_{h} u = u H(u^{\alpha},x,s).
\end{equation}
%
%
\end{rem}

\begin{lem}
\label{lem:hessian}
 Let $\mathcal{H}_{h}$ be a differential operator on
$C^\infty_{>0}(N)$ defined by
    \begin{equation}
     \mathcal{H}_{h} v =\sum r_{i}\frac{H_{h}^{v^{a_{i}}}  }
{v^{a_{i}}},
        \label{eq:hessian.1}
    \end{equation}
    $\zeta :=\sum r_{i}a_{i}$ and
    $\eta:= \sum r_{i}a_{i}^{2}$, where the indices extend from
    $1$ to $l \in \N $ and any $r_{i},a_{i} \in \mathbb{R}$. Hence,
  \begin{equation}
  \label{eq:hessian.2}
  \mathcal{H}_{h} v  = (\eta-\zeta)\frac{1}{v^{2}} dv \otimes dv  +
\zeta \frac{1}{v}
  H_{h}^{v},
  \end{equation}
  where $\otimes $ is the usual tensorial product.
  If furthermore, $\zeta \neq 0$ and $\eta \neq 0$, then
    \begin{equation}
     \mathcal{H}_{h} v= \beta \frac{H_{h}^{v^{\frac{1}{\alpha}}}
}
     {v^{\frac{1}{\alpha}}},
        \label{eq:hessian.3}
    \end{equation}
    where $\alpha = \displaystyle \frac{\zeta}{\eta}$ and
    $\beta = \displaystyle \frac{\zeta^{2}}{\eta}$.

\end{lem}

\section{The problem \textbf{(Pb-sc)} - Existence of solutions}
\label{sec:The problem (Pb-sc)}

\emph{Throughout this section, we will assume that B is not only a
Riemannian manifold of dimension $m \ge 3$, but also ``compact"
and connected.  We further assume that $F$ is a pseudo-Riemannian
manifold of dimension $k \ge 0 $ with constant scalar curvature
$S_F \le 0$. Moreover, we will assume that $\mu \neq \mu_{sc}$.}

\smallskip

\noindent Hence,  we will concentrate our attention on the relations
\eqref{eq:sc principal 1}, \eqref{eq:alpha principal} and
\eqref{eq:beta principal} by applying Theorem \ref{thm:scurv conf
warped m ge 2} (i).

\medskip

Let $\lambda_{1}$ denote the principal eigenvalue of the operator
\begin{equation}
\label{eq:lambda_1 (S_B)}
   L(\cdot) = -\beta \Delta_{B}(\cdot) + S_{B}
   (\cdot),
\end{equation}
and $u_{1} \in C^{\infty}_{>0}(B)$ be the corresponding positive
eigenfunction with $\|u_1\|_\infty = 1$, where $\beta $ is as in
Theorem \ref{thm:scurv conf warped m ge 2}.

\medskip

First of all, we will state some results about uniqueness and
non-existence of positive solutions for Equation \eqref{eq:scalar
curv bcwp1} under the latter hypothesis.

\smallskip

About the former, we adapt Lemma 3.3 in \cite[p.
525]{Ambrosetti-Brezis-Cerami94} to our situation (for a detailed
proof see \cite{Ambrosetti-Brezis-Cerami94}, \cite[Method II, p.
103]{Brezis-Kamin92} and also \cite{Shi-Yao05}).

\begin{lem}
\label{lem:lemma 3.3 ABC} Let $f \in C^{0}(\mathbb{R}_{>0})$ such
that $t^{-1}f(t)$ is decreasing. If $v$ and $w$ satisfy
\begin{equation}
  \begin{array}{c}
    -\beta \Delta_{B}v + S_{B} v  \le f(v),\\
  v \in C^{\infty}_{>0}(B),
  \end{array}
 \label{eq:lemma 3.3 ABC 1}
\end{equation}
and
\begin{equation}
  \begin{array}{c}
    -\beta \Delta_{B}w + S_{B} w  \ge f(w),\\
  w \in C^{\infty}_{>0}(B),
  \end{array}
 \label{eq:lemma 3.3 ABC 2}
\end{equation}
then $w \ge v$ on $B$.
\end{lem}

\begin{proof}
Let $\theta(t)$ be a smooth nondecreasing function such that
$\theta(t)\equiv 0$ for $t \le 0$ and $\theta(t)\equiv 1$ for $t
\ge 1$. Thus for all $\epsilon >0$,
\begin{equation*}
   \theta_\epsilon(t):=\theta\left(\frac{t}{\epsilon}\right)
\end{equation*}
is smooth, nondecreasing, nonnegative and $\theta(t)\equiv 0$ for $t
\le 0$ and $\theta(t)\equiv 1$ for $t \ge \epsilon$. Furthermore
$\gamma_\epsilon(t):=\int_0^t s\theta_\epsilon^\prime(s) ds$
satisfies $0 \le \gamma_\epsilon(t) \le \epsilon,$ for any $t \in
\mathbb{R}$.

\noindent On the other hand, since $(B,g_B)$ is a compact Riemannian
manifold without boundary and $\beta>0$, like in \cite[Lemma 3.3, p.
526]{Ambrosetti-Brezis-Cerami94} there results
\begin{equation}\label{eq:lemma 3.3 ABC 3}
  \int_B [-v \beta\Delta_B w + w \beta\Delta_B v] \theta_\epsilon(v-w) dv_{g_B}\le
  \int_B [-\beta\Delta_B v]\gamma_\epsilon(v-w)dv_{g_B}.
\end{equation}
Hence, by the above considerations about $\theta_\epsilon$ and
$\gamma_\epsilon$, \eqref{eq:lemma 3.3 ABC 3} implies that
\begin{equation}\label{eq:lemma 3.3 ABC 4}
  \int_B [-v \beta\Delta_B w + w \beta\Delta_B v] \theta_\epsilon(v-w) dv_{g_B}\le
  \epsilon \int_{[-\beta\Delta_B v \ge 0]} [-\beta\Delta_B v ]  dv_{g_B}.
\end{equation}

\noindent Now, by applying \eqref{eq:lemma 3.3 ABC 1} and
\eqref{eq:lemma 3.3 ABC 2} there results
\begin{equation}\label{eq:lemma 3.3 ABC 5}
  -v \beta\Delta_B w + w \beta\Delta_B v = v L w - w L v \ge v f (w) - w f
  (v)=vw\left[\frac{f(w)}{w}-\frac{f(v)}{v}\right].
\end{equation}

\noindent Thus by combining \eqref{eq:lemma 3.3 ABC 5} and
\eqref{eq:lemma 3.3 ABC 4}, as $\epsilon \rightarrow 0^+$ we led
to
\begin{equation}\label{eq:lemma 3.3 ABC 6}
  \int_{[v>w]} vw\left[\frac{f(w)}{w}-\frac{f(v)}{v}\right] dv_{g_B}\le 0
\end{equation}
and conclude the proof like in \cite[Lemma 3.3, p.
526-527]{Ambrosetti-Brezis-Cerami94}. But $\displaystyle
\frac{f(v)}{v}<\frac{f(w)}{w}$ on $[v>w]$ and hence $\meas[v>w]=0$;
thus $v \le w$. \footnote{$\meas$ denotes the usual $g_B-$measure on
the compact Riemannian manifold $(B_m,g_B)$}
\end{proof}

\begin{cor}
\label{cor:corollary-lemma 3.3 ABC} Let $f \in
C^{0}(\mathbb{R}_{>0})$ such that $t^{-1}f(t)$ is decreasing. Then
\begin{equation}
  \begin{array}{c}
    -\beta \Delta_{B}v + S_{B} v  = f(v),\\
  v \in C^{\infty}_{>0}(B)
  \end{array}
 \label{eq:corollary-lemma 3.3 ABC 1}
\end{equation}
has at most one solution.
\end{cor}

\begin{proof}
Assume that $v$ and $w$ are two solutions of
\eqref{eq:corollary-lemma 3.3 ABC 1}. Then by applying Lemma
\ref{lem:lemma 3.3 ABC} firstly with $v$ and $w$, and conversely
with $w$ and $v$, the conclusion is proved.
\end{proof}

\begin{rem}\label{rem:singular}
Notice that Lemma \ref{lem:lemma 3.3 ABC} and Corollary
\ref{cor:corollary-lemma 3.3 ABC} allow the function $f \in
C^{0}(\mathbb{R}_{>0})$ to be singular at $0$.
\end{rem}

\smallskip

Related to the non-existence of smooth positive solutions for
Equation \eqref{eq:scalar curv bcwp1}, we will state an easy result
under the general hypothesis of this section.

\begin{prop}\label{prop:non existence}
If either $\,\, \max_B S_B \,\, \le \,\, \inf_{u \in
\mathbb{R}_{>0}} u^{2\mu\alpha}(\lambda - S_{F}u^{-2\alpha}) \,$
or \\ $\min_B S_B \ge \sup_{u \in \mathbb{R}_{>0}}
u^{2\mu\alpha}(\lambda  - S_{F}u^{-2\alpha})$,
%
%
then \eqref{eq:scalar curv bcwp1} has no solution in
$C^\infty_{>0} (B)$.
\end{prop}

\begin{proof}It is sufficient to apply the
maximum principle with some easy adjustments to the particular
involved coefficients.
\end{proof}

\medskip

\emph{$\bullet$ The case of scalar flat fiber, i.e. $S_F=0$.}

\medskip

In this case, the term containing the nonlinearity $u^{2(\mu -
1)\alpha + 1}$ becomes non-influent in \eqref{eq:scalar curv bcwp1},
thus \textbf{(Pb-sc)} equivalently results to the study of existence
of solutions for the problem:
\begin{equation}
  \begin{array}{c}
    -\beta \Delta_{B}u + S_{B} u  = \lambda u^{2 \mu \alpha + 1},\\
  u \in C^{\infty}_{>0}(B),
  \end{array}
 \label{eq:sc principal SF=0}
\end{equation}
where $\lambda$ is a real parameter (i.e., it is the searched
constant scalar curvature) and $\psi=u^\alpha$.

\begin{rem}\label{rem:difference of homogeneity}
\footnote{Along this article we consider the sign function defined
by $\sign = \chi_{(0,+\infty)} - \chi_{(-\infty,0)}$, where $\chi_A$
is the characteristic function of the set $A$. } Let $p \in
\mathbb{R} \setminus \{1\}$ and $(\lambda_0,u_0) \in
(\mathbb{R}\setminus \{0\}) \times C^{\infty}_{>0}(B)$ be a solution
of
\begin{equation}
  \begin{array}{c}
    -\beta \Delta_{B}u + S_{B} u  = \lambda u^{p},\\
  u \in C^{\infty}_{>0}(B).
  \end{array}
 \label{eq:sc principal SF=0 p}
\end{equation}
Hence, by the difference of homogeneity between both members of
\eqref{eq:sc principal SF=0}, it is easy to show that if $\lambda
\in \mathbb{R}$ satisfies $\sign (\lambda) = \sign (\lambda_0)$,
then $(\lambda,u_\lambda )$ is a solution of \eqref{eq:sc
principal SF=0 p}, where $u_\lambda=t_\lambda u_0$ and $t_\lambda
= \displaystyle
\left(\frac{\lambda}{\lambda_0}\right)^{\frac{1}{1-p}} $.

\noindent Thus by \eqref{eq:sc principal SF=0}, we obtain
geometrically: if the parameter $\mu $ is given in a way that $p:=2
\mu \alpha +1 \neq 1$ and \twistpar{B}{F}{\psi_0^\mu}{\psi_0} has
constant scalar curvature $\lambda_0 \neq 0$, then for any $\lambda
\in \mathbb{R}$ verifying $\sign (\lambda) = \sign (\lambda_0)$,
there results that \twistpar{B}{F}{\psi_\lambda^\mu}{\psi_\lambda}
is of scalar curvature $\lambda$, where $\psi_\lambda =
t_\lambda^\alpha \psi_0$ and $t_\lambda$ given as above.
\end{rem}

\begin{thm}
\label{thm:mu=0 linear} $(Case \colon \,\, \mu=0)$
The scalar curvature of a \emph{\bcwpar{\psi,0}} of base $B$ and
fiber $F$ (i.e., a singly warped product $B \times_\psi F$) is a
constant $\lambda$ if and only if $ \lambda=\lambda_{1}$ and
$\psi$ is a positive multiple of $u_1^{\frac{2}{k+1}}$ (i.e.,
$\psi = t u_1^{\frac{2}{k+1}}$ for some $t \in \mathbb{R}_{>0}$).
\end{thm}

\begin{proof}
First of all note that $\mu =0 $ implies $\displaystyle
\alpha=\frac{2}{k+1}$. On the other hand, in this case, the problem
\eqref{eq:sc principal SF=0} is linear, so it is sufficient to apply
the well known results about the principal eigenvalue and its
associated eigenfunctions of operators like \eqref{eq:lambda_1
(S_B)} in a suitable setting.
\end{proof}

\begin{thm}
\label{thm:musc < mu < 0 - sub-linear}
$(Case \colon \,\, {\mu_{sc} < \mu < 0})$
The scalar curvature of a \emph{\bcwpar{\psi,\mu}} of base $B$ and
fiber $F$ is a constant $\lambda$, only if $\sign (\lambda) = \sign
(\lambda_{1})$. Furthermore,
\begin{enumerate}
\item if $\lambda = 0$ then there exists $\psi \in C^{\infty}_{>0}(B)$
such that \twistpar{B}{F}{\psi^\mu}{\psi} has constant scalar
curvature $0$ if and only if $\lambda_1=0$. Moreover, such
$\psi$'s are the positive multiples of $u_1^\alpha$, i.e.
$tu_1^\alpha$, $t \in \mathbb{R}_{>0}$.
\item if $\lambda > 0$ then there exists $\psi \in C^{\infty}_{>0}(B)$
such that \twistpar{B}{F}{\psi^\mu}{\psi} has constant scalar
curvature $\lambda$ if and only if $\lambda_1>0$. In this case, the
solution $\psi$ is unique.
\item if $\lambda < 0$ then there exists $\psi \in C^{\infty}_{>0}(B)$
such that \twistpar{B}{F}{\psi^\mu}{\psi} has constant scalar
curvature $\lambda$ when $\lambda_1<0$ and is close enough to $0$.
\end{enumerate}
\end{thm}

\begin{proof} The condition $\mu_{sc} < \mu
< 0$ implies that $0<p:=2 \mu \alpha + 1<1$, i.e., the problem
\eqref{eq:sc principal SF=0} is sublinear. Thus, to prove the
theorem one can use variational arguments as in
\cite{CotiZelatiDobarroMusina97} (alternatively, degree theoretic
arguments as in \cite{Ambrosetti-Hess80} or bifurcation theory as in
\cite{DobarroLamiDozo87}).

\noindent We observe that in order to obtain the positivity of the
solutions required in \eqref{eq:sc principal SF=0}, one may apply
the maximum principle for the case of $\lambda > 0$ and the
antimaximum principle for the case of $\lambda < 0$.

\noindent The uniqueness for $\lambda > 0$ is a consequence of
Corollary \ref{cor:corollary-lemma 3.3 ABC}.
\end{proof}

\begin{rem}\label{rem:Hebey}
In order to consider the next case we introduce the following
notation. For a given $p$ such that $1<p\le p_Y$, let
\begin{equation}\label{eq:kappa p}
        \kappa_p := \inf_{v \in \mathcal{H}_p}
        \int_B \left(|\nabla^B v|^2 + \frac{S_B}{\beta}
        v^2\right)dv_{g_B},
\end{equation}
    where
\begin{equation*}
    \mathcal{H}_p := \left\{ v \in H^1(B) : \int_B |v|^{p+1}dv_{g_B} =1
    \right\}.
\end{equation*}
Now, we consider the following two cases.

\begin{description}
    \item[$(1<p < p_Y)$]
In this case by adapting \cite[Theorem 1.3]{Hebey03}, there exists
$u_p \in C_{>0}^\infty(B)$ such that $(\beta\kappa_p,u_p)$ is a
solution of \eqref{eq:sc principal SF=0 p}
%
%
and $\displaystyle \int_B u_p^{p+1}dv_{g_B} =1$.
%
%
    \item[$(p= p_Y)$]For this specific and important value, analogously
    to \cite[\S2]{Hebey03}, we distinguish
three subcases along the study of our problem \eqref{eq:sc
principal SF=0 p}, in correspondence with the $\sign
(\kappa_{p_Y})$.
\begin{description}
    \item[$\underline{\kappa_{p_Y}=0}$] in this case, there exists
    $u_{p_Y} \in C_{>0}^\infty(B)$
such that $(0,u_{p_Y})$ is a solution of \eqref{eq:sc principal
SF=0 p} and $\displaystyle \int_B u_{p_Y}^{p_Y +1}dv_{g_B} =1$.
    \item[$\underline{\kappa_{p_Y}<0}$] here there exists
    $u_{p_Y} \in C_{>0}^\infty(B)$
such that $(\beta\kappa_{p_Y},u_{p_Y})$ is a solution of
\eqref{eq:sc principal SF=0 p} and $\displaystyle \int_B
u_{p_Y}^{p_Y +1}dv_{g_B} =1$.
    \item[$\underline{\kappa_{p_Y}>0}$] this is a more difficult
    case, let $K_m$ be the sharp Euclidean Sobolev constant
\begin{equation}\label{eq:sharp constant}
    K_m=\sqrt{\displaystyle
    \frac{4}{m(m-2)\omega_m^{\frac{2}{m}}}},
\end{equation}
where $\omega_m$ is the volume of the unit $m-$sphere. Thus, if
\begin{equation}
\label{eq:sharp condition_1} \kappa_{p_Y}<\frac{1}{K_m^2},
\end{equation}
then there exists $u_{p_Y} \in C_{>0}^\infty(B)$ such that
$(\beta\kappa_{p_Y},u_{p_Y})$ is a solution of \eqref{eq:sc
principal SF=0 p} and $\displaystyle \int_B u_{p_Y}^{p_Y +1}dv_{g_B}
=1$. Furthermore, the condition
\begin{equation}
\label{eq:sharp condition_2} \kappa_{p_Y}\le \frac{1}{K_m^2}
\end{equation}
is sharp by \cite{Hebey03}, so that this is independent of the
underlying manifold and the potential considered.

The equality case in \eqref{eq:sharp condition_2} is discussed in
\cite{Hebey-Vaugon01}.
\end{description}
\end{description}
This results allow to establish the following two theorems.
\end{rem}

\begin{thm}
\label{thm:mupY < mu < musc - (super-linear,sub-critical)} $(Cases
\colon \,\, {\mu_{p_Y}< \mu < \mu_{sc}\textrm{ or }0 < \mu })$
There exists $\psi \in C_{>0}^\infty (B)$ such that the scalar
curvature of \twistpar{B}{F}{\psi^\mu}{\psi} is a constant $\lambda$
if and only if $\sign (\lambda) = \sign (\kappa_p)$ where $p:=2 \mu
\alpha + 1$ and $\kappa_p$ is given by \eqref{eq:kappa p}.
Furthermore if $\lambda < 0,$ then the solution $\psi$ is unique.
\end{thm}

\begin{proof}
The conditions $(\mu_{p_Y}< \mu < \mu_{sc}\textrm{ or }0 < \mu )$
imply that $1<p:=2 \mu \alpha + 1<p_Y$, i.e. the problem
\eqref{eq:sc principal SF=0} is superlinear but subcritical with
respect to the Sobolev immersion theorem
(see \cite[Remark 5.5]{DobarroUnal04-2}).
%
By recalling that $\psi = u^\alpha$, it is sufficient to prove
that follows.

Let $u_p$ be defined as in the case of $(1<p<p_Y)$ in Remark
\ref{rem:Hebey}. If $(\lambda,u)$ is a solution of \eqref{eq:sc
principal SF=0}, then multiplying \eqref{eq:sc principal SF=0} by
$u_{p}$ and integrating by parts there results
\begin{equation}\label{}
  \beta\kappa_p \int_{B} u_{p} u dv_{g_B}= \lambda \int_{B} u_{p} u^{p} dv_{g_B}.
\end{equation}
Thus $\sign (\lambda) = \sign (\kappa_p)$ since $\beta$, $u_{p}$ and
$u$ are all positive.

Conversely, if $\lambda $ is a real constant such that $\sign
(\lambda) = \sign (\kappa_p)\neq 0$, then by \textit{Remark
\ref{rem:difference of homogeneity}}, $(\lambda,u_\lambda )$ is a
solution of \eqref{eq:sc principal SF=0}, where
$u_\lambda=t_\lambda u_p$ and $t_\lambda = \displaystyle
\left(\frac{\lambda}{\beta \kappa_p}\right)^{\frac{1}{1-p}} $.

On the other side, if $\lambda = \kappa_p =0$, then $(0,u_p)$ is a
solution of \eqref{eq:sc principal SF=0}.

\noindent Since $1<p$, the uniqueness for $\lambda < 0$ is a
consequence of Corollary \ref{cor:corollary-lemma 3.3 ABC}.
\end{proof}

\begin{thm}
\label{thm:mupY  - critical} $(Cases \colon \,\, \mu = \mu_{p_Y})$
If there exists $\psi \in C_{>0}^\infty (B)$ such that the scalar
curvature of \twistpar{B}{F}{\psi^{\mu_{p_Y}}}{\psi} is a constant
$\lambda$, then $\sign (\lambda) = \sign (\kappa_{p_Y})$.
Furthermore, if $\lambda \in \mathbb{R}$ verifying $\sign (\lambda)
= \sign (\kappa_{p_Y})$ and \eqref{eq:sharp condition_1}, then there
exists $\psi \in C_{>0}^\infty (B)$ such that the scalar curvature
of \twistpar{B}{F}{\psi^{\mu_{p_Y}}}{\psi} is $\lambda$. Besides, if
$\lambda \in \mathbb{R}$ is negative, then there exists at most one
$\psi \in C_{>0}^\infty (B)$ such that the scalar curvature of
\twistpar{B}{F}{\psi^{\mu_{p_Y}}}{\psi} is $\lambda$.
\end{thm}

\begin{proof}
The proof is similar to that of Theorem \ref{thm:mupY < mu < musc -
(super-linear,sub-critical)}, but follows from the application of
the case of $(p=p_Y)$ in Remark \ref{rem:Hebey}. Like above, the
uniqueness of $\lambda < 0$ is a consequence of Corollary
\ref{cor:corollary-lemma 3.3 ABC}.
\end{proof}

\medskip

In the next proposition including the supercritical case, we will
apply the following result (see also \cite[p.99]{Taylor96}).

\begin{lem}\label{lem: taylor 3p99}
Let $(N_n,g_N)$ be a compact connected Riemannian manifold without
boundary of dimension $n \ge 2 $ and $\Delta_{g_N}$ be the
corresponding Laplace-Beltrami operator. Consider the equation of
the form
\begin{equation}\label{eq:taylor 3p99}
    \begin{array}{c}
    -\Delta_{g_N}u=f(\cdot,u),\\
    u \in C^\infty_{>0}(N)
    \end{array}
\end{equation}
where $f \in C^\infty(N \times \mathbb{R}_{>0})$. If there exist
$a_0$ and $a_1 \in \mathbb{R}_{>0}$ such that
\begin{equation}\label{eq:taylor condition 3p99}
    \begin{array}{c}
     u <a_0 \Rightarrow f(\cdot,u) >0\\
{\text and} \\
  u >a_1 \Rightarrow f(\cdot,u) <0,
    \end{array}
\end{equation}
then \eqref{eq:taylor 3p99} has a solution satisfying $a_0 \le u
\le a_1$.
\end{lem}

\medskip

\begin{prop}\label{prop:supercritical SB<0 SF=0}
$(Cases \colon \,\, {-\infty< \mu < \mu_{sc}\textrm{ or }0 < \mu })$
If $\max S_B < 0$, then for all $\lambda < 0$ there exists $\psi \in
C_{>0}^\infty (B)$ such that the scalar curvature of
\twistpar{B}{F}{\psi^\mu}{\psi} is the constant $\lambda$.
Furthermore, the solution $\psi$ is unique.
\end{prop}

\begin{proof}
The conditions $(-\infty < \mu < \mu_{sc}\textrm{ or }0 < \mu )$
imply that $1<p:=2 \mu \alpha + 1$.

\noindent On the other hand, since $B$ is compact, by taking
\begin{equation*}\label{}
    f(.,u)=
-S_{B}(\cdot) u + \lambda u^{p}=(-S_B + \lambda u^{p-1})u,
\end{equation*}
we obtain that $\lim _{u \longrightarrow 0^+}f(\cdot,u)= 0^+$ and
$\lim _{u \longrightarrow +\infty}f(\cdot,u)= -\infty$. Thus
\eqref{eq:taylor condition 3p99} is verified.

\noindent Hence, the proposition is proved by applying Lemma
\ref{lem: taylor 3p99} on $(B_m,g_B)$. Notice that $a_0$ can take
positive values and eventually gets close enough to $0^+$ due to the
condition of $\lim _{u \longrightarrow 0^+}f(\cdot,u)$, and
consequently the corresponding solution results positive.

\noindent Again, since $\lambda < 0$ and $1<p$ the uniqueness is a
consequence of Corollary \ref{cor:corollary-lemma 3.3 ABC}.
\end{proof}

\begin{proof}(of Theorem \ref{thm:Pb-sc SF0}) This is an immediate
consequence of the above results.
\end{proof}

\emph{$\bullet$ The case of a fiber with negative constant scalar
curvature, i.e. $S_F<0$.}

\medskip

Here, the \textbf{(Pb-sc)} becomes equivalent to the study of the
existence for the problem
\begin{equation}
  \begin{array}{c}
    -\beta \Delta_{B}u + S_{B} u  = \lambda u^{p} - S_{F}u^{q},\\
  u \in C^{\infty}_{>0}(B),
  \end{array}
 \label{eq:sc principal SF < 0}
\end{equation}
where $\lambda$ is a real parameter (i.e., the searched constant
scalar curvature), $\psi=u^\alpha$, $p=2 \mu \alpha + 1$ and
$q=2(\mu - 1)\alpha + 1$.

\begin{rem} \label{rem:lambda position} Let $u$ be a solution of
\eqref{eq:sc principal SF < 0}.
\begin{itemize}
    \item[{\bf (i)}] If $\lambda_{1} \leq 0$,
  then $\lambda < 0$. Indeed, multiplying the equation in
  \eqref{eq:sc principal SF < 0} by $u_{1}$
  and integrating by parts there results:
\begin{equation}\label{}
  \lambda_{1} \int_{B} u_{1} u dv_{g_B}+ S_{F} \int_{B} u_{1} u^{q} dv_{g_B}=
  \lambda \int_{B} u_{1} u^{p}dv_{g_B},
\end{equation}
where $u_{1}$ and $u$ are positive.

  \item[{\bf (ii)}] If $\lambda = 0$, then $
  \lambda_{1} >0$.

  \item[{\bf (iii)}] If $\mu = 0$ (the warped product case), then
  $\lambda < \lambda_1$.
  These cases have been studied in
  \cite{DobarroLamiDozo87,{CotiZelatiDobarroMusina97}}.

 \item[{\bf (iv)}] If $\mu = 1$ (the Yamabe problem for the usual
  product with conformal factor in $C^{\infty}_{>0}(B)$),
  there results $\sign (\lambda) = \sign (\lambda_{1} +
  S_{F})$.
\end{itemize}
\end{rem}

An immediate consequence of Remark \ref{rem:lambda position} is the
following lemma.

\begin{lem}\label{rem:non existence SF<0}
Let $B$ and $F$ be given like in Theorem \ref{thm:scurv conf warped
m ge 2}(i). Suppose further that $B$ is a compact connected
Riemannian manifold and $F$ is a pseudo-Riemannian manifold of
constant scalar curvature $S_F<0$. If $\lambda \ge 0$ and $\lambda_1
\le 0$ (for instance when $S_B \le 0$ on $B$), then there is no
$\psi \in C^\infty _{>0}(B)$ such that the scalar curvature of
\twistpar{B}{F}{\psi^{\mu}}{\psi} is $\lambda$.
\end{lem}

\begin{thm}
\label{thm:concave-convex} \cite[Rows 6 and 8 in Table
4]{DobarroUnal04-2} Under the hypothesis of Theorem \ref{thm:scurv
conf warped m ge 2}(i), let $B$ be a compact connected Riemannian
manifold and $F$ be a pseudo-Riemannian manifold of constant scalar
curvature $S_F<0$. Suppose that ``$(m,k) \in D$ and $\mu \in (0,1)
$" or ``$(m,k) \in \mathcal{C}D$ and $\mu \in (0,1) \cap
\mathcal{C}[\mu_{-} , \mu_{+}]$".
%
%
\begin{enumerate}
    \item If $\lambda_1 \le 0$, then $\lambda \in \mathbb{R}$ is
    the scalar curvature of a \twistpar{B}{F}{\psi^{\mu}}{\psi}
    if and only if $\lambda < 0$.

    \item If $\lambda_1 > 0$, then there exists $\overline{\Lambda}
    \in \mathbb{R}_{>0}$ such that $\lambda \in \mathbb{R}\setminus
    \{\overline{\Lambda}\}$ is the scalar curvature of a
    \twistpar{B}{F}{\psi^{\mu}}{\psi} if and only if
    $\lambda < \overline{\Lambda}$.
%
%
\end{enumerate}
Furthermore if $\lambda \le 0$, then there exists at most one $\psi
\in C^{\infty}_{>0}(B)$ such that \twistpar{B}{F}{\psi^{\mu}}{\psi}
has scalar curvature $\lambda$.
\end{thm}

\begin{proof}
The proof of this theorem is the subject matter of \S 5.
\end{proof}

\medskip

Once again we make use of Lemma \ref{lem: taylor 3p99} for the
next theorem about the singular case and the following
propositions.

\begin{thm}
\label{thm:p superlin q sing} \cite[Row 7 Table 4]{DobarroUnal04-2}
Under the hypothesis of Theorem \ref{thm:scurv conf warped m ge
2}(i), let $B$ be a compact connected Riemannian manifold and $F$ be
a pseudo-Riemannian manifold of constant scalar curvature $S_F<0$.
Suppose that ``$(m,k) \in \mathcal{C}D$ and $\mu \in (0,1) \cap \;
(\mu_{-} , \mu_{+})$",
%
%
then for any $\lambda <0$ there exists $\psi \in C^\infty _{>0}(B)$
such that the scalar curvature of \twistpar{B}{F}{\psi^{\mu}}{\psi}
is $\lambda$. Furthermore the solution $\psi $ is unique.
\end{thm}

\begin{proof}
First of all note that the conditions ``$(m,k) \in \mathcal{C}D$
and $\mu \in (0,1) \cap \; (\mu_{-} , \mu_{+})$" imply that $q <
0$ and $1<p$ , i.e. the problem \eqref{eq:sc principal SF < 0} is
superlinear in $p$ but singular in $q$.

On the other hand, since $B$ is compact, taking
\begin{equation*}\label{}
    f(.,u)=
-S_{B}(\cdot) u + \lambda u^{p} - S_{F}u^{q}=[(-S_B(\cdot) +
\lambda u^{p-1})u^{1-q}-S_F]u^q,
\end{equation*}
there result $\lim _{u \longrightarrow 0^+}f(\cdot,u)= +\infty$
and $\lim _{u \longrightarrow +\infty}f(\cdot,u)= -\infty$. Thus
\eqref{eq:taylor condition 3p99} is verified.

Thus by an application of Lemma \ref{lem: taylor 3p99} for
$(B_m,g_B),$ we conclude the proof for the existence part.

The uniqueness part just follows from Corollary
\ref{cor:corollary-lemma 3.3 ABC}.
\end{proof}

\begin{rem}\label{rem:taylor for concave-convex}
We observe that the arguments applied in the proof of Theorem
\ref{thm:p superlin q sing} can be adjusted to the case of a
compact connected Riemannian manifold $B$ with $0 \le q < 1 < p$,
$\lambda < 0$ and $S_F<0$, so that some of the situations included
in Theorem \ref{thm:concave-convex}. However, both argumentations
are compatible but different.
\end{rem}

\begin{proof}(of Theorem \ref{thm:Pb-sc SF<0}) This is an immediate
consequence of the above results.
\end{proof}

The approach in the next propositions is similar to Proposition
\ref{prop:supercritical SB<0 SF=0} and Theorem \ref{thm:p superlin q
sing}.

\begin{prop}\label{prop:supercritical SB<0 SF<0}
%
\cite[Row 10 Table 4]{DobarroUnal04-2}
Let $1< \mu <+\infty$. If $\max S_B < 0$, then for all $\lambda < 0$
there exists $\psi \in C_{>0}^\infty (B)$ such that the scalar
curvature of \twistpar{B}{F}{\psi^\mu}{\psi} is the constant
$\lambda$.
%
%
%
\end{prop}

\begin{proof}
The condition $1 < \mu < +\infty$ implies that $1<q<p$.

On the other hand, since $B$ is compact, taking
\begin{equation*}\label{}
    f(.,u)=
-S_{B}(\cdot) u + \lambda u^{p} - S_{F}u^{q}=[-S_B(\cdot) +
(\lambda u^{p-q}-S_F)u^{q-1}]u,
\end{equation*}
there result $\lim _{u \longrightarrow 0^+}f(\cdot,u)= 0^+$ and
$\lim _{u \longrightarrow +\infty}f(\cdot,u)= -\infty$. Thus
\eqref{eq:taylor condition 3p99} is satisfied.

Thus an elementary application of Lemma \ref{lem: taylor 3p99} for
$(B_m,g_B)$ proves the proposition.
\end{proof}

\begin{prop}\label{prop:pqsublin or sublin-sing SB>0 SF<0}
\cite[Rows 2, 4 and 3 in Table 4]{DobarroUnal04-2}
Let either ``$(m,k) \in D$ and $\mu \in (\mu_{sc},0) $" or
``$(m,k) \in \mathcal{C}D$ and $\mu \in (\mu_{sc},0) \cap
\mathcal{C}[\mu_{-}, \mu_{+}]$" or ``$(m,k) \in \mathcal{C}D$ and
$\mu \in (\mu_{sc},0) \cap (\mu_{-} , \mu_{+})$". If $\min S_B >
0$, then for all $\lambda \le 0$ there exists a smooth function
$\psi \in C_{>0}^\infty (B)$ such that the scalar curvature of
\twistpar{B}{F}{\psi^\mu}{\psi} is the constant $\lambda$.
%
%
%
\end{prop}

\begin{proof}
If either ``$(m,k) \in D$ and $\mu \in (\mu_{sc},0) $" or ``$(m,k)
\in \mathcal{C}D$ and $\mu \in (\mu_{sc},0) \cap
\mathcal{C}[\mu_{-}, \mu_{+}]$", then $0<q<p<1$.

On the other hand, since $B$ is compact, taking
\begin{equation*}\label{}
    f(.,u)=
-S_{B}(\cdot) u + \lambda u^{p} - S_{F}u^{q}=[-S_B(\cdot)u^{1-q} +
\lambda u^{p-q}-S_F]u^q,
\end{equation*}
there result $\lim _{u \longrightarrow 0^+}f(\cdot,u)= 0^+$ and
$\lim _{u \longrightarrow +\infty}f(\cdot,u)= -\infty$. Thus
\eqref{eq:taylor condition 3p99} is verified and again we can apply
Lemma \ref{lem: taylor 3p99} for $(B_m,g_B)$.

If ``$(m,k) \in \mathcal{C}D$ and $\mu \in (\mu_{sc},0) \cap
(\mu_{-} , \mu_{+})$", then $q<0<p<1$. Considering the limits as
above, $\lim _{u \longrightarrow 0^+}f(\cdot,u)= +\infty$ and $\lim
_{u \longrightarrow +\infty}f(\cdot,u)= -\infty$. So, an application
of Lemma \ref{lem: taylor 3p99} concludes the proof.
\end{proof}

\begin{rem}\label{rem:sign of SB}
Notice that in Theorems \ref{thm:concave-convex} and \ref{thm:p
superlin q sing} we do not assume hypothesis related to the sign of
$S_B(\cdot)$, unlike in Propositions \ref{prop:supercritical SB<0
SF=0}, \ref{prop:supercritical SB<0 SF<0} and \ref{prop:pqsublin or
sublin-sing SB>0 SF<0}.
\end{rem}

\begin{prop}\label{prop:non homog cases}
\cite[Rows 5 and 9 in Table 4]{DobarroUnal04-2}
Let $(m,k) \in \mathcal{C}D$ be.
\begin{enumerate}
    \item If either ``$\mu \in
\left(-\displaystyle\frac{k}{m-1},0\right)\cap \{\mu_{-} ,
\mu_{+}\}$ and $\min S_B >0$" or ``$\mu \in (0,1)\cap \{\mu_{-} ,
\mu_{+}\}$", then for all $\lambda < 0$ there exists a smooth
function $\psi \in C_{>0}^\infty (B)$ such that the scalar curvature
of \twistpar{B}{F}{\psi^\mu}{\psi} is the constant $\lambda$. In the
second case, $\psi$ is also unique .

    \item If either ``$\mu \in
\left(-\displaystyle\frac{k}{m-1},0\right)\cap \{\mu_{-} ,
\mu_{+}\}$" or ``$\mu \in (0,1)\cap \{\mu_{-} , \mu_{+}\}$" and
furthermore $\lambda_1 >0$, then there exists a smooth function
$\psi \in C_{>0}^\infty (B)$ such that the scalar curvature of
\twistpar{B}{F}{\psi^\mu}{\psi} is $0$.
\end{enumerate}
\end{prop}

\begin{proof}
In both cases $q=0$, so by considering
    \begin{equation*}\label{}
    f(.,u)=
    -S_{B}(\cdot) u + \lambda u^{p} - S_{F},
    \end{equation*}
the proof of (1) follows as in the latter propositions, while that
of (2) is a consequence of the linear theory and the maximum
principle.
\end{proof}

\begin{rem}\label{rem: surface fiber DL}
Finally, we observe a particular result about the cases studied in
\cite{DobarroLamiDozo87}. If $\mu = 0$, then $p=1$ and
$\displaystyle q=1-2\alpha =\frac{k-3}{k+1}$. When the dimension of
the fiber is $k=2$, the exponent $\displaystyle q=-\frac{1}{3}$. So,
writing the involved equation as
    \begin{equation*}\label{}
    -\frac{8}{3}\Delta_B u = f(.,u)=
    -S_{B}(\cdot) u + \lambda u - S_{F}u^{-\frac{1}{3}}
    \end{equation*}
and by applying Lemma \ref{lem: taylor 3p99} as above, we obtain
that if $\lambda < \min S_B$, then there exists a smooth function
$\psi \in C_{>0}^\infty (B)$ such that the scalar curvature of $B
\times_\psi F$ is the constant $\lambda$. Furthermore, by Corollary
\ref{cor:corollary-lemma 3.3 ABC} such $\psi$ is unique (see
\cite{DobarroLamiDozo87, CotiZelatiDobarroMusina97} and
\cite{CrandallRabinowitzTartar76}).
\end{rem}


\section{Proof of the Theorem \ref{thm:concave-convex}}
\label{sec:concave-convex}

The subject matter of this section is the proof of the Theorem
\ref{thm:concave-convex}, so we naturally assume its hypothesis.

\noindent Most of the time, we need to specify the dependence of
$\lambda $ of \eqref{eq:sc principal SF < 0}, we will do that by
writing $\eqref{eq:sc principal SF < 0}_\lambda$. \noindent
Furthermore, we will denote the right hand side of $\eqref{eq:sc
principal SF < 0}_\lambda$  by $f_\lambda (t):=\lambda t^{p} -
S_{F}t^{q}$.

\smallskip

The conditions either ``$(m,k) \in D$ and $\mu \in (0,1) $" or
``$(m,k) \in \mathcal{C}D$ and $\mu \in (0,1) \cap
\mathcal{C}[\mu_{-} , \mu_{+}]$", imply that $0<q<1<p$. But the
type of nonlinearity in the right hand side of $\eqref{eq:sc
principal SF < 0}_\lambda$ changes with the $\sign \, \lambda$,
i.e. it is purely concave for $\lambda <0$ and concave-convex for
$\lambda
>0$.

\noindent The uniqueness for $\lambda \le 0$ is again a
consequence of Corollary \ref{cor:corollary-lemma 3.3 ABC}.

\noindent In order to prove the existence of a solution for
$\eqref{eq:sc principal SF < 0}_\lambda$ with $\sign \,\lambda
\neq 0$, we adapt the approach of sub and upper solutions in
\cite{Ambrosetti-Brezis-Cerami94}.
\medskip

\noindent Thus, the proof of Theorem \ref{thm:concave-convex} will
be an immediate  consequence of the results that follows.

\medskip

\begin{lem}
\label{lem:lambda=0, 0<q<1} $\eqref{eq:sc principal SF < 0}_0$ has
a solution if and only if $\lambda_1>0$.
\end{lem}

\begin{proof}
This situation is included in the results of the second case of
Theorem \ref{thm:musc < mu < 0 - sub-linear} by replacing $-S_F$
with $\lambda$ (see \cite[Proposition
3.1]{CotiZelatiDobarroMusina97}).
\end{proof}

\begin{lem}
\label{lem:Lambda range $S_F < 0, 0<q<1<p$}
 Let us assume that
$\{\lambda : \eqref{eq:sc principal SF < 0}_\lambda \textrm{ has a
solution}\}$ is non-empty and define
\begin{equation}\label{}
    \overline{\Lambda}=\sup \{\lambda : \eqref{eq:sc principal SF < 0}_\lambda \textrm{ has a solution}\}.
\end{equation}

\begin{itemize}
    \item[{\bf (i)}] If $\lambda_1 \le 0$, then  $\overline{\Lambda} \le 0$.
    \item[{\bf (ii)}] If $\lambda_1 > 0$, then there exists $\overline{\lambda} >0$ finite such
that $\overline{\Lambda} \le \overline{\lambda}$.
\end{itemize}
\end{lem}

\begin{proof}
$ $
\begin{itemize}
    \item[{\bf (i)}] It is sufficient to observe \textit{Remark \ref{rem:lambda position}
    i}.
    \item[{\bf (ii)}] Like in \cite{Ambrosetti-Brezis-Cerami94}, let
$\overline{\lambda}>0$ such that
\begin{equation}\label{eq:$f;S_F < 0, 0<q<1<p$}
    \lambda_1 t < \overline{\lambda} t^p - S_F t^q, \forall t \in
    \mathbb{R}, t>0.
\end{equation}
Thus, if $(\lambda,u)$ is a solution of $\eqref{eq:sc principal SF <
0}_\lambda$, then

\begin{equation*}\label{eq:}
  \lambda \int_{B} u_{1} u^{p} - S_{F} \int_{B} u_{1} u^{q}=\int_{B} \lambda_{1} u_{1} u
  < \overline{\lambda} \int_{B} u_{1} u^{p} - S_{F} \int_{B} u_{1} u^{q},
\end{equation*}
so $\lambda < \overline{\lambda}$.
\end{itemize}
\end{proof}

\begin{lem}
\label{lem:existence $S_F < 0, 0<q<1<p, lambda > 0$}
 Let
\begin{equation}\label{eq:overlineLambda}
    \overline{\Lambda}=\sup \{\lambda : \eqref{eq:sc principal SF < 0}_\lambda
    \textrm{ has a solution}\}.
\end{equation}

\begin{itemize}
    \item[{\bf (i)}] Let $E \in \mathbb{R}_{>0}$.
    There exist $0<\lambda_0=\lambda_0(E)$  and $0<M=M(E,\lambda_0)$
    such that $\forall \lambda: 0 < \lambda \le \lambda_0$,
    so we have
\begin{equation}\label{eq:supersol1}
    0<E\frac{f_\lambda (EM)}{EM}<1.
\end{equation}

    \item[{\bf (ii)}] If $\lambda_1 > 0$, then  $\{\lambda >0 :
    \eqref{eq:sc principal SF < 0}_\lambda
    \textrm{ has a solution}\} \neq \emptyset$.
    As a consequence of that, $\overline{\Lambda}$ is finite.

    \item[{\bf (iii)}] If $\lambda_1 > 0$, then for all
    $0<\lambda<\overline{\Lambda}$ there exists
    a solution of the problem $\eqref{eq:sc principal SF < 0}_\lambda$.
\end{itemize}
\end{lem}

  \begin{figure}
   \epsfig{file=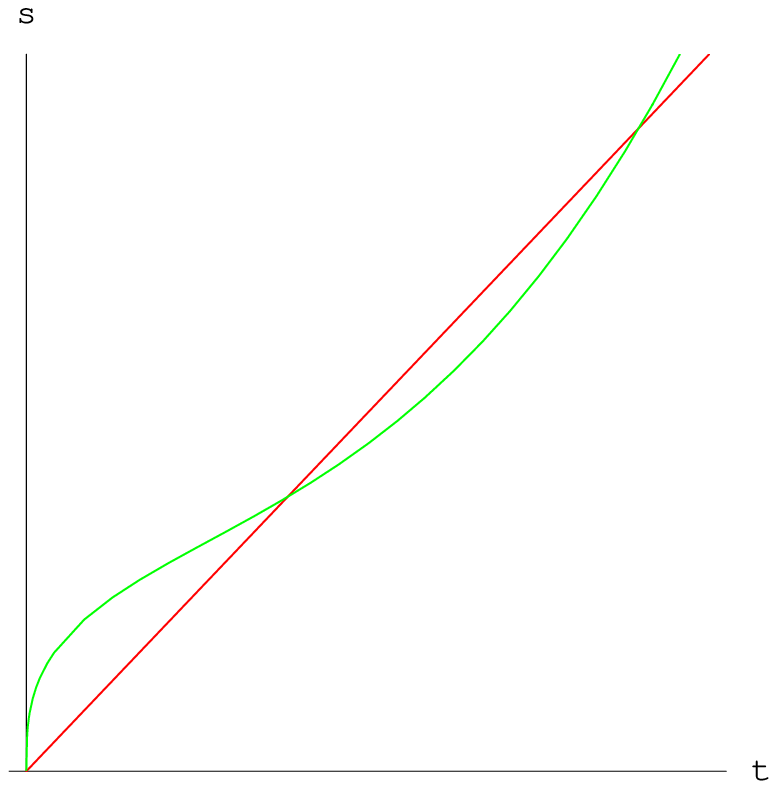,width=7.2cm,angle=0}
   \begin{small}
   \caption{
   The nonlinearity $f_\lambda$ in Lemma
   \ref{lem:existence $S_F < 0, 0<q<1<p, lambda > 0$},
   i.e.  $0<q<1<p, S_F < 0, \lambda_1 > 0, \lambda > 0$.
   }
   \end{small}
  \end{figure}

\begin{proof}$ $
\begin{itemize}
    \item[{\bf (i)}] For any $0<\lambda<\lambda_0$
    \begin{align*}
0<g_\lambda(r):=E\frac{f_\lambda (Er)}{Er}
&= E r^{q-1} ( \lambda E^{p-1}r^{p-q} - S_F E^{q-1})\\
&<E r^{q-1} ( \lambda_0 E^{p-1}r^{p-q} - S_F E^{q-1}).
    \end{align*}
    It is easy to
    see that
    \begin{equation*}
r_0=\left(\frac{S_F}{\lambda_0}\frac{q-1}{p-1}\right)^{\frac{1}{p-q}}\frac{1}{E}
    \end{equation*}
    is a minimum point for $g_{\lambda_0}$ and
    \begin{equation*}
g_{\lambda_0} (r_0)=
E\left(\frac{S_F}{\lambda_0}\frac{q-1}{p-1}\right)^{\frac{q-1}{p-q}}S_F
\left[\frac{q-1}{p-1}-1\right]\rightarrow 0^+, \textrm{ as }
\lambda_0 \rightarrow 0^+.
    \end{equation*}
    Hence there exist $0<\lambda_0=\lambda_0(E)$ and $0<M=M(E,\lambda_0)$ such
    that \eqref{eq:supersol1} is verified.
    \item[{\bf (ii)}] Since $\lambda_1 >0$, by the maximum principle,
    there exists a solution $e \in C^{\infty}_{>0}(B)$ of
    \begin{equation}\label{eq:unit sol 1} L_B(e) = -\beta \Delta_{B} e + S_{B} e =1.
    \end{equation}
    Then, applying item (i) above with $E=\|e\|_\infty$
    there exists $0<\lambda_0=\lambda_0(\|e\|_\infty)$
    and $0<M=M(\|e\|_\infty,\lambda_0)$ such
    that $\forall \lambda$ with $0<\lambda\le\lambda_0$
    we have that
    \begin{equation}\label{eq:supersol2}
     L_B (Me) = M \ge f_\lambda(Me),
    \end{equation}
    hence $Me$ is a supersolution of $\eqref{eq:sc principal SF < 0}_\lambda$.

On the other hand, since $\check{u}_1:=\inf u_1>0$, for all
$\lambda
>0$
    \begin{equation}\label{eq:supersol3}
     \epsilon^{-1}f_\lambda(\epsilon \check{u}_1)=
     \epsilon^{q-1} [\lambda \epsilon^{p-q} \check{u}_1^p-S_F \check{u}_1^q]\rightarrow
+\infty, \textrm{ as } \epsilon \rightarrow 0^+.
    \end{equation}
Furthermore, note that $f_\lambda$ is nondecreasing when $\lambda
>0$. Hence for any $0<\lambda$ there exists a small enough $0<\epsilon$ verifying
    \begin{equation}\label{eq:supersol4}
     L_B (\epsilon u_1) = \epsilon \lambda_1 u_1 \le \epsilon \lambda_1 \|u_1\|_\infty
     \le f_\lambda(\epsilon \check{u}_1)\le f_\lambda(\epsilon u_1),
    \end{equation}
thus $\epsilon u_1$ is a subsolution of $\eqref{eq:sc principal SF
< 0}_\lambda$.

Then for any $0<\lambda<\lambda_0$, (taking eventually $0<\epsilon
$ smaller if necessary), we have that the above constructed couple
sub super solution satisfies
    \begin{equation}\label{eq:supersol5}
     \epsilon u_1 < M e.
    \end{equation}
Now, by applying the monotone iteration scheme, we have that
\\
$\{\lambda >0 : \eqref{eq:sc principal SF < 0}_\lambda \textrm{ has
a solution}\} \neq \emptyset$. Furthermore by Lemma \ref{lem:Lambda
range $S_F < 0, 0<q<1<p$} (ii) there results $\overline{\Lambda}$ is
finite.

   \item[{\bf (iii)}] The proof of this item is completely analogous to
    \textit{Lemma 3.2} in \cite{Ambrosetti-Brezis-Cerami94}. We
    will rewrite this to be self contained.

    Given $\lambda < \overline{\Lambda}$, let $u_\nu$ be a
    solution of $\eqref{eq:sc principal SF < 0}_\nu$ with
    $\lambda<\nu<\overline{\Lambda}$. Then $u_\nu$ is a
    supersolution of $\eqref{eq:sc principal SF < 0}_\lambda$
    and for small enough $0< \epsilon $, the subsolution $\epsilon u_1$ of
    $\eqref{eq:sc principal SF < 0}_\lambda$ verifies $\epsilon u_1 < u_\nu$,
    then as above $\eqref{eq:sc principal SF < 0}_\lambda$ has a solution.

\end{itemize}
\end{proof}

\begin{lem}
\label{lem: norma infty $S_F < 0, 0<q<1<p, lambda < 0$} For any
$\lambda<0$, there exists $\gamma_\lambda > 0$ such that
$\|u\|_\infty \le \gamma_\lambda$ for any solution $u$ of
$\eqref{eq:sc principal SF < 0}_\lambda$. Furthermore if $S_B$ is
nonnegative, then positive zero of $f_\lambda$ can be choose as
$\gamma_\lambda$.
\end{lem}

\begin{proof}
Define $\check{S}_B:=\min S_B$ (recall that $B$ is compact). There
are two different situations, namely.
    \begin{itemize}
        \item $\underline{0 \le \check{S}_B }$:
        since there
        exists $x_1 \in B$ such that $u(x_1)=\|u\|_\infty$ and $0 \le -\beta \Delta_{B}u(x_1) =
        - S_{B}(x_1) \|u\|_\infty + \lambda \|u\|_\infty^{p}-
        S_{F}\|u\|_\infty^{q}$, there results $\|u\|_\infty \le
        \gamma_\lambda$, where $\gamma_\lambda$ is the strictly positive zero
        of $f_\lambda$.
        \item $\underline{\check{S}_B < 0}$: we consider
        $\widetilde{f}_\lambda(t):=\lambda t^p - S_F t^q - \check{S}_B
        t$. Now, our problem $\eqref{eq:sc principal SF < 0}_\lambda$ is equivalent to
\begin{equation*}
\label{eq:constant scalar curvature m ge 2:2}
  \begin{array}{c}
    -\beta \Delta_{B}u + (S_{B}-\check{S}_B) u  = \widetilde{f}_\lambda(u),\\
  u \in C^{\infty}_{>0}(B).
  \end{array}
\end{equation*}
But here the potential of $(S_{B}-\check{S}_B)$ is non negative and
the function $\widetilde{f}_\lambda$ has the same behavior of
$f_\lambda$ with a positive zero $\widetilde{\gamma}_\lambda $ on
the right side of the positive zero $\gamma_\lambda $ of
$f_\lambda$. Thus, repeating the argument for the case of $
\check{S}_B \ge 0$, we proved $\|u\|_\infty \le
        \tilde{\gamma}_\lambda$.
    \end{itemize}
\end{proof}

\begin{lem}
\label{lem:existence $S_F < 0, 0<q<1<p, lambda le 0$}
 Let $\lambda_1 > 0$. Then
for all $\lambda < 0$ there exists
    a solution of $\eqref{eq:sc principal SF < 0}_\lambda$.
\end{lem}

\begin{proof}   We will apply again the
monotone iteration scheme. Define $\check{S}_B:=\min S_B$ (note
that $B$ is compact).

  \begin{figure}
    \label{fig:Lemma }
   \epsfig{file=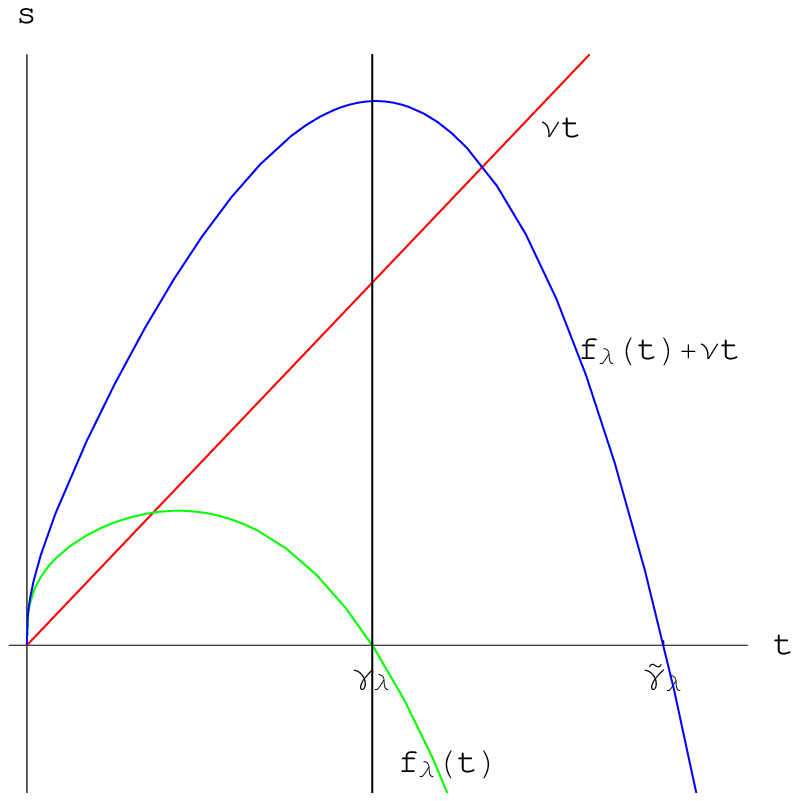,width=7.2cm,angle=0}
    \begin{small}
   \caption{The nonlinearity in Lemma \ref{lem:existence $S_F < 0, 0<q<1<p, lambda le 0$}
   ,
   i.e.  $0<q<1<p$, $S_F < 0$, $\lambda_1 > 0$, $\lambda < 0$.
   }
    \end{small}
  \end{figure}

    \begin{itemize}
    \item $\underline{0 \le \check{S}_B }$:
Clearly, the strictly positive zero $\gamma_\lambda$ of
$f_\lambda$ is a supersolution of
\begin{equation}
\label{eq:constant scalar curvature m ge 2 + nu}
    -\beta \Delta_{B}u + (S_{B}+\nu) u  = f_\lambda (u)  + \nu u,
\end{equation}
for all $\nu \in \mathbb{R}$.

On the other hand, for $0<\epsilon = \epsilon(\lambda)$ small
enough,
    \begin{equation}\label{eq:subsolution lambda<0 1}
L_B (\epsilon u_1) = \epsilon \lambda_1 u_1 \le f_\lambda(\epsilon
u_1).
    \end{equation}
Then $\epsilon u_1$ is a subsolution of \eqref{eq:constant scalar
curvature m ge 2 + nu} for all $\nu \in \mathbb{R}$.

By taking $\varepsilon$ possibly smaller, we also have
\begin{equation}\label{eq:well order}
   0< \epsilon u_1 < \gamma_\lambda.
\end{equation}

We note that for large enough values of $\nu \in \mathbb{R}_{>0}$,
the nonlinearity on the right hand side of \eqref{eq:constant
scalar curvature m ge 2 + nu}, namely $f_\lambda (t)  + \nu t$, is
an increasing function on $[0, \gamma_\lambda]$.

Thus applying the monotone iteration scheme we obtain a strictly
positive solution of \eqref{eq:constant scalar curvature m ge 2 +
nu}, and hence a solution of $\eqref{eq:sc principal SF <
0}_\lambda$ (see \cite{Amann72}, \cite{Amann76},
\cite{PLLions82}).

%
%
%

    \item $\underline{\check{S}_B < 0}$: In this case, like in
    Lemma \ref{lem: norma infty $S_F < 0, 0<q<1<p, lambda < 0$} we
    consider $\widetilde{f}_\lambda(t):=\lambda t^p - S_F t^q - \check{S}_B
        t$. Then, the problem $\eqref{eq:sc principal SF < 0}_\lambda$ is equivalent to
\begin{equation}
\label{eq:constant scalar curvature m ge 2:3}
  \begin{array}{c}
    -\beta \Delta_{B}u + (S_{B}-\check{S}_B) u  = \widetilde{f}_\lambda(u),\\
  u \in C^{\infty}_{>0}(B),
  \end{array}
\end{equation}
where the potential is nonnegative and the function
$\widetilde{f}_\lambda$ has a similar behavior to $f_\lambda$ with a
positive zero $\widetilde{\gamma}_\lambda $ on the right side of the
positive zero $\gamma_\lambda $ of $f_\lambda$.

Here, it is clear that $\tilde{\gamma}_\lambda$ is a positive
supersolution of
\begin{equation}
\label{eq:constant scalar curvature m ge 2:3 + nu}
    -\beta \Delta_{B}u + (S_{B}-\check{S}_B +\nu) u  =
    \widetilde{f}_\lambda(u) + \nu u,
\end{equation}
for all $\nu \in \mathbb{R}$. Hence, we complete the proof similarly
to the case of $\check{S}_B \ge 0$.
    \end{itemize}
\end{proof}

\begin{lem}
\label{lem: non existence $S_F < 0, 0<q<1<p, lambda_1 le 0,lambda le
0$} Let $\lambda_1 \le 0$, $\lambda < 0$, $\check{S}_B:=\min S_B$
and also let $\gamma_\lambda$ be a positive zero of $f_\lambda $ and
$\tilde{\gamma}_\lambda$ be a positive zero of
$\tilde{f}_\lambda:=f_\lambda - \check{S}_B id_{\mathbb{R}_{\ge
0}}$. Then there exists a solution $u$ of $\eqref{eq:sc principal SF
< 0}_\lambda$. Furthermore any solution of $\eqref{eq:sc principal
SF < 0}_\lambda$ satisfies $\gamma_\lambda \le \|u\|_\infty \le
\tilde{\gamma}_\lambda$.
\end{lem}

\begin{proof}
First of all we observe that if $S_B \equiv 0$ (so $\lambda_1 =
0$), then $u \equiv \gamma_\lambda$ is the searched solution of
$\eqref{eq:sc principal SF < 0}_\lambda$.

Now, we assume that $S_B \not\equiv 0$. Since $\lambda_1 \le 0$,
there results $\check{S}_B < 0$. In this case, one can notice that
$0<\gamma_\lambda<\tilde{\gamma}_\lambda$.

\noindent On the other hand, the problem $\eqref{eq:sc principal
SF < 0}_\lambda$ is equivalent to
\begin{equation}
\label{eq:constant scalar curvature m ge 2:4}
  \begin{array}{c}
    -\beta \Delta_{B}u + (S_{B}-\check{S}_B) u  = \widetilde{f}_\lambda(u),\\
  u \in C^{\infty}_{>0}(B).
  \end{array}
\end{equation}

\noindent By the second part of the proof of Lemma \ref{lem: norma
infty $S_F < 0, 0<q<1<p, lambda < 0$}, if $u$ is a solution of
$\eqref{eq:sc principal SF < 0}_\lambda$ (or equivalently
\eqref{eq:constant scalar curvature m ge 2:4}), then $\|u\|_\infty
\le \tilde{\gamma}_\lambda $. Besides, since
\begin{equation*}\label{eq:}
  \int_{B} u_{1} (f_\lambda \circ u) =\lambda_{1} \int_{B}  u_{1} u,
\end{equation*}
$u, u_1 >0$ and $\lambda_{1} \le 0$ results $\gamma_\lambda \le
\|u\|_\infty $.

\noindent From this point on, the proof of the existence of
solutions for \eqref{eq:constant scalar curvature m ge 2:4} follows
the lines of the second part of Lemma \ref{lem:existence $S_F < 0,
0<q<1<p, lambda le 0$}.
\end{proof}

\section{Conclusions and future directions}
\label{sec:Conclusions}

Now, we would like to summarize the content of the paper and to
propose our future plans on this topic.

We remark to the reader that several computations and proofs,
along with other complementary results mentioned in this article
and references can be obtained in \cite{DobarroUnal04-2}. We have
chosen this procedure to avoid the involved long computations.

In brief, we introduced and studied curvature properties of a
particular family of warped products of two pseudo-Riemannian
manifolds which we called as a \textit{base conformal warped
product}. Roughly speaking the metric of such a product is a mixture
of a conformal metric on the base and a warped metric. We
concentrated our attention on a special subclass of this structure,
where there is a specific relation between the \emph{conformal
factor} $c$ and the \emph{warping function} $w$, namely $c=w^\mu$
with $\mu $ a real parameter.

\noindent As we mentioned in \S 1 and the first part of \S 2, these
kinds of metrics and considerations about their curvatures are very
frequent in different physical areas, for instance theory of general
relativity, extra-dimension theories (Kaluza-Klein,
Randall-Sundrum), string and super-gravity theories; also in global
analysis for example in the study of the spectrum of
Laplace-Beltrami operators on $p$-forms, etc.

More precisely, in Theorems \ref{thm:global riki} and
\ref{thm:sca-c}, we obtained the classical relations among the
different involved Ricci tensors (respectively, scalar curvatures)
for metrics of the form $c^{2}g_B \oplus w^2 g_F$. Then the study
of particular families of either scalar or tensorial nonlinear
partial differential operators on pseudo-Riemannian manifolds (see
Lemmas \ref{m-lem} and \ref{lem:hessian}) allowed us to find
reduced expressions of the Ricci tensor and scalar curvature for
metrics as above with $c=w^\mu$, where $\mu $ a real parameter
(see Theorems \ref{thm:ricci sbcwp m ge 3} and \ref{thm:scurv conf
warped m ge 2}). The operated reductions can be considered as
generalizations of those used by Yamabe in \cite{Yamabe63} in
order to obtain the transformation law of the scalar curvature
under a conformal change in the metric and those used in
\cite{DobarroLamiDozo87} with the aim to obtain a suitable
relation among the involved scalar curvatures in a singly warped
product (see also \cite{Lelong-Ferrand76} for other particular
application and our study on multiply warped products in
\cite{DobarroUnal04}).

In \S4 and 5, under the hypothesis that $(B,g_B)$ be a ``compact"
and connected Riemannian manifold of dimension $m \ge 3$ and
$(F,g_F)$ be a pseudo - Riemannian manifold of dimension $k \ge 0 $
with constant scalar curvature $S_F$,  we dealt with the problem
\textbf{(Pb-sc)}. This question leads us to analyze the existence
and uniqueness of solutions for nonlinear elliptic partial
differential equations with several kinds of nonlinearities. The
type of nonlinearity changes with the value of the real parameter
$\mu$ and the sign of $S_F$. In this article, we concentrated our
attention to the cases of constant scalar curvature $S_F \le 0$ and
accordingly the central results are Theorems \ref{thm:Pb-sc SF0} and
\ref{thm:Pb-sc SF<0}. Although our results are partial so that there
are more cases to study in forthcoming works, we obtained also other
complementary results under more restricted hypothesis about the
sign of the scalar curvature of the base.

Throughout our study, we meet several types of partial
differential equations. Among them, most important ones are those
with concave-convex nonlinearities and the one so called
Lichnerowicz-York equation. About the former, we deal with the
existence of solutions and leave the question of multiplicity of
solutions to a forthcoming study.

We observe that the previous problems as well as the study of the
Einstein equation on \textit{base conformal warped products},
\emph{\bcwpar{\psi,\mu}'s} and their \textit{generalizations to
multi-fiber cases}, give rise to a reach family of interesting
problems in differential geometry and physics (see for instance, the
several recent works of R. Argurio, J. P. Gauntlett, M. O. Katanaev,
H. Kodama, J. Maldacena, H. -J. Schmidt, A. Strominger, K. Uzawa, P.
S. Wesson among many others) and in nonlinear analysis (see the
different works of A. Ambrosetti, T. Aubin, I. Choquet-Bruat, J.
Escobar, E. Hebey, J. Isenberg, A. Malchiodi, D. Pollack, R. Schoen,
S. -T. Yau among others).

\appendix{}
\numberwithin{equation}{section}
\renewcommand{\theequation}{\Alph{section}.\arabic{equation}}
\section{ }
\label{appendix:A}

Let us assume the hypothesis of Theorem \ref{thm:scurv conf warped m
ge 2} (i), the dimensions of the base $m \ge 2$ and of the fiber $k
\ge 1$. In order to describe the classification of the type of
nonlinearities involved in \eqref{eq:sc principal 1}, we will
introduce some notation (for a complete study of these
nonlinearities see \cite[Section 5]{DobarroUnal04-2}). The example
in Figure 1 will help the reader to clarify the notation.

\medskip

Note that the denominator in \eqref{eq:alpha principal} is
\begin{equation}\label{eq:scalar eta}
 \eta := (m-1)(m-2)\mu^{2} + 2(m-2)k\mu +(k+1)k
\end{equation}
and verifies $\eta > 0 \textrm{ for all } \mu \in \mathbb{R}$. Thus
$\alpha $ in \eqref{eq:alpha principal} is positive if and only if
$\displaystyle \mu > -\frac{k}{m-1}$ and by the hypothesis
$\displaystyle \mu \neq -\frac{k}{m-1}$ in Theorem \ref{thm:scurv
conf warped m ge 2} (i), results $\alpha \neq 0$.

\medskip

We now introduce the following notation:
\begin{equation}\label{eq:p,q}
  \begin{array}{ll}
    p= p(m,k,\mu)= & 2 \mu \alpha + 1 \textrm{ and }\\
    q= q(m,k,\mu)= & 2(\mu - 1)\alpha + 1 = p - 2 \alpha,
  \end{array}
\end{equation}
where $\alpha$ is defined by \eqref{eq:alpha principal}.

\medskip

Thus, for all $ m, k, \mu $ given as above, $p$ is positive. Indeed,
by \eqref{eq:scalar eta}, $p > 0$ if and only if $\varpi > 0$, where
\begin{equation*}\label{}
 \begin{array}{ll}
   \varpi &:= \varpi (m,k,\mu) \\
      &:=4 \mu [k + (m-1)\mu] + (m-1)(m-2)\mu^{2} + 2(m-2)k\mu
+(k+1)k \\
          &=(m-1)(m+2)\mu^{2} + 2 m k \mu  + (k+1)k.
 \end{array}
\end{equation*}
But $\discr{\varpi} \le -4km^{2} \le -16$ and $m>1$, so $\varpi
>0.$

\medskip

Unlike $p$, $q$ changes sign depending on $m$ and $k$.
Furthermore, it is important to determine the position of $p$ and
$q$ with respect to $1$ as a function of $m$ and $k$. In order to
do that, we define
\begin{equation}\label{eq:m,k condition}
  D:=\{(m,k)\in \mathbb{N}_{\ge 2}\times\mathbb{N}_{\ge 1}:
\discr{\varrho(m,k,\cdot)}<0\},
\end{equation}
where $\mathbb{N}_{\ge l} := \{j \in \mathbb{N}: j \ge l\}$ and
\begin{equation*}\label{}
 \begin{array}{ll}
   \varrho &:= \varrho (m,k,\mu) \\
       &:=4 (\mu-1) [k + (m-1)\mu] + (m-1)(m-2)\mu^{2} +
2(m-2)k\mu +(k+1)k \\
           &=(m-1)(m+2)\mu^{2} + 2 (m k -2(m-1))\mu  + (k-3)k.
 \end{array}
\end{equation*}
Note that by \eqref{eq:scalar eta}, $q>0$ if and only if $\varrho
>0$. Furthermore $q=0$ if and only if $\varrho
=0$. But here $\discr{\varrho(m,k,\cdot)}$ changes its sign as a
function of $m$ and $k$.


We adopt here the notation in \cite[\textsc{Table
4}]{DobarroUnal04-2} below, namely $\mathcal{C}D = (\mathbb{N}_{\ge
2} \times \mathbb{N}_{\ge 1}) \setminus D$ if $D \subseteq
\mathbb{N}_{\ge 2} \times \mathbb{N}_{\ge 1}$ and $\mathcal{C}I =
\mathbb{R} \setminus I$ if $I \subseteq \mathbb{R}$. Thus, if $(m,k)
\in \mathcal{C}D$, let $\mu_{-}$ and $\mu_{+}$ two roots (eventually
one, see \cite[Remark 5.3]{DobarroUnal04-2}) of $q$, $\mu_{-} \le
\mu_{+}$. Besides, if $\discr{\varrho(m,k,\cdot)}>0 $, then
$\mu_{-}<0$; whereas $\mu_{+}$ can take any sign.

\begin{figure}
    \label{fig:notation}
   \epsfig{file=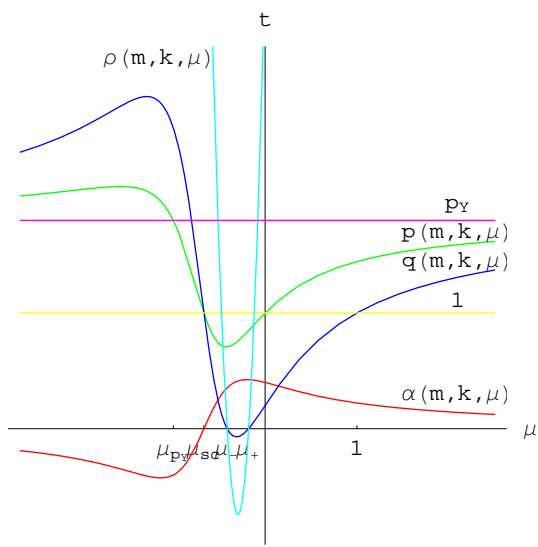,width=7.2cm,angle=0}
   \begin{small}
   \caption{Example:
   $(m,k)=(7,4) \in \mathcal{C}D$
   }
   \end{small}
\end{figure}

\providecommand{\bysame}{\leavevmode\hbox
to3em{\hrulefill}\thinspace}

\end{document}